\documentclass[hidelinks,onefignum,onetabnum]{siamart250211}

\usepackage{lipsum}
\usepackage{amsfonts}
\usepackage{graphicx,xcolor}
\usepackage{epstopdf}
\ifpdf
  \DeclareGraphicsExtensions{.eps,.pdf,.png,.jpg}
\else
  \DeclareGraphicsExtensions{.eps}
\fi
 
\usepackage[margin=1in]{geometry}

\headers{Cross$^2$-DEIM and Tucker-AA }{D. Appel\"{o} and Y. Cheng}

\title{A new cross approximation for Tucker tensors and its application in Tucker-Anderson Acceleration}

\author{Daniel Appel\"{o} \thanks{Department of Mathematics, Virginia Tech, Blacksburg, VA 24061 U.S.A. (\textsf{appelo@vt.edu}). Research supported by U.S. Department of Energy, Office of Science, Advanced Scientific Computing Research (ASCR), under Award Number DE-SC0025424, NSF DMS-2208164, NSF DMS-2436319, and Virginia Tech.} 
\and Yingda Cheng \thanks{Department of Mathematics, Virginia Tech, Blacksburg, VA 24061 U.S.A. (\textsf{yingda@vt.edu}). Research is supported by DOE grant DE-SC0023164, AFOSR grant FA9550-25-1-0154 and Virginia Tech. This material is based upon work supported by the National Science Foundation under Grant No.
DMS-1928930 while the authors were in residence at the Simons Laufer Mathematical Sciences
Institute in Berkeley, California, during the Fall 2025 semester.} }

\usepackage{comment}

\usepackage[compatible]{algpseudocode} 
\renewcommand{\algorithmiccomment}[1]{\bgroup\hfill\small\#~#1\egroup}

\graphicspath{{figures/}}

\begin{document}
\maketitle

\begin{abstract} 
This paper proposes two new algorithms related to the Tucker tensor format. The first  method is a new cross approximation for Tucker tensors, which we call Cross$^2$-DEIM. Cross$^2$-DEIM is an iterative method that uses a fiber sampling strategy, sampling $O(r)$ fibers in each mode, where $r$ denotes the target rank. The fibers are selected based on the discrete empirical interpolation method (DEIM). Cross$^2$-DEIM resemblances the Fiber Sampling Tucker Decomposition (FSTD)2 approximation, and has favorable computational scaling compared to existing methods in the literature. We demonstrate good performance of Cross$^2$-DEIM in terms of iteration count and intermediate memory. First we design a fast direct Poisson solver based on Cross$^2$-DEIM and the fast Fourier transform. This solver can be used as a stand alone or as a preconditioner for low-rank solvers for elliptic problems.

The second method is a low-rank solver for nonlinear tensor equation in Tucker format by Anderson acceleration (AA), which we call Tucker-AA. Tucker-AA   is an extension of low-rank AA (lrAA) proposed in our prior work \cite{appelo2025lraa} for low-rank solution to nonlinear matrix equation. We apply Cross$^2$-DEIM with warm-start in Tucker-AA to deal with the nonlinearity in the equation. We apply low-rank operations in AA, and by an appropriate rank truncation strategy, we are able to control the intermediate rank growth. We demonstrated the   performance for Tucker-AA for approximate solutions nonlinear PDEs in 3D.
\end{abstract}



\section{Introduction}
In this paper, we consider numerical algorithms for Tucker tensors. The Tucker decomposition \cite{hitchcock1927expression,tucker1966some} is a widely used tensor format that decomposes a tensor into a core tensor and factor matrices. It also serves as a building block for more sophisticated tensor formats such as tree tensor networks. Tucker tensors have been shown to be effective and robust in compressing scientific data. For a $d$-dimensional tensor with $n$ elements in each dimension, a Tucker tensor stores $O(r^d+dnr)$ degrees of freedom where $r$ denotes the multi-linear rank. If the rank $r$ is low this results in significant saving compared to storing full tensor which requires $O(n^d)$ storage.

This paper is motivated by the need to solve nonlinear tensor equations of the kind
\begin{equation}
\label{eq:nle}
H(X)=X, \quad X \in \mathbb{R}^{n_1 \times n_2 \ldots \times n_d},
\end{equation}
when the tensor $X$ is represented using the Tucker tensor format. Equation \eqref{eq:nle} models tensor equations arising from discretizations of differential equations in moderate dimensions, where Tucker tensor can offer good compression ratio. This is an extension of our prior work \cite{appelo2025lraa}: lrAA, low-rank Anderson Acceleration for nonlinear matrix equations, to the Tucker tensor case. In \eqref{eq:nle}, the function $H(\cdot)$ is a element-wise nonlinear function, for example, it could come from a discretization of a nonlinear PDE. 

 Recently, there have been great interests in computing high dimensional differential equations using low-rank tensor based approaches \cite{hackbusch2012tensor,khoromskij2018tensor,bachmayr2023low}. Broadly speaking, there are three types of low-rank methods for solving \eqref{eq:nle}, those based on iteration combined with low-rank truncation, those based on optimization and those based on greedy rank enrichment. We refer the readers to \cite{grasedyck2013literature,bachmayr2023low} for details. For the problem that we are interested in, solving linear and nonlinear systems in Tucker format, we mention the work in \cite{kressner2016preconditioned}, where a Riemannian optimization approach  with preconditioning is proposed for Tucker and tensor train manifolds. Further, the work in \cite{georgieva2019greedy} proposed a greedy method for solving linear system in Tucker format, the contribution \cite{luo2023low} developed a Riemannian Gauss-Newton method for Tucker tensor estimation. We also mention the recent work \cite{nakao2025reduced} where an implicit solver for 3D linear tensor equations in Tucker format, motivated by the basis update and Galerkin method \cite{ceruti2022unconventional} for time-dependent linear advection-diffusion equation, is developed.   

The goal of this paper is to design an iterative method with low-rank truncation. We chose the framework of Anderson acceleration  (AA) \cite{anderson1965iterative, saad2024acceleration} because it is amenable for controlling intermediate rank during the iterations. Compared to Newton method, AA does not need access to the Jacobian, which may not be low-rank. Further, AA accommodates a natural warm-start strategy for   cross approximation, which is the costly operation that resolves the nonlinearity.
  In \cite{appelo2025lraa}, we proposed a novel framework for computing adaptive low-rank solution to nonlinear matrix equations by  lrAA.   Given a prescribed tolerance, lrAA outputs numerical solution in its singular value decomposition (SVD) form based on a tolerance criteria. This paper aims to address the extension of lrAA to the Tucker tensor format.  

Because cross approximation is a bottleneck for computational efficiency, to extend the ideas to Tucker tensor format in high dimensions, in this paper, we develop a new  efficient cross approximation for Tucker tensor. We call it Cross$^2$-DEIM. Compared to  existing cross approximations for Tucker tensor in the literature \cite{ahmadi2021cross}, it has a favorable scaling in terms of computational cost  and storage. It is an iterative method which uses a fiber sampling strategy guided by discrete empirical interpolation method (DEIM)  index selection.  Specifically, by using DEIM on approximate factor matrices, we determine the main index sets.  We devise a fast strategy by cross approximation for the unfoldings of a sub-sampled tensor to decide the complement index sets. The name Cross$^2$-DEIM comes from the fact that  two cross approximations are used to choose the main and complement index guided by DEIM. For a $d$-dimensional tensor with length $n$  in each mode, we show that the computational cost of each iteration of Cross$^2$-DEIM to be $O(dn r^2+r^{d+1})$ with $O(dnr+r^d)$ evaluations of the tensor entries, where $r$ denotes the intermediate rank. We show the iteration numbers of Cross$^2$-DEIM stay low  in practice. We use the Cross$^2$-DEIM algorithm in Tucker-Anderson Acceleration (Tucker-AA) and observe good numerical performance on benchmark problems arising from numerical approximations to nonlinear PDEs.

The rest of the paper is organized as follows. In Section \ref{sec:tuckercross}, we present a review of mathematical notations for Tucker tensors and introduce Cross$^2$-DEIM approximation. Section \ref{sec:tuckeraa} presents Tucker-AA method. Section \ref{sec:num} contains numerical experiments of the two methods. In Section \ref{sec:conclude} we conclude the paper.

\section{A new cross approximation for Tucker tensor}
\label{sec:tuckercross}
In this section, we review standard tensor notations and introduce Cross$^2$-DEIM: a new cross approximation for Tucker tensors.

\subsection{Tucker tensor format}
\label{sec:tucker}
We follow the standard notation for Tucker tensor \cite{kolda2009tensor}, which is briefly reviewed below. We consider a $d$-dimensional tensor $X \in \mathbb{R}^{n_1\times n_2 \times \cdots n_d}.$ To annotate the possible elements in the $i$th dimension we define $\mathcal{N}_i=\{1, \ldots, n_i \}, i=1, \ldots d$.  

The  $i$-th mode matricization of $X$ is defined by
$$
X_{(i)} \in  \mathbb{R}^{n_i\times \Pi_{j \ne i} n_j}.
$$
This is simply the unfolding of $X$ into a matrix, where  the row index corresponds to the $i$-th mode, and the column index corresponds to the linear lexicographically ordered index over  the other modes. We denote the rank of $X_{(i)}$ to be $r_i$ and the multilinear rank of the tensor $X$ to be $\mathbf{r}=(r_1, \ldots, r_d)$. We define the tensor norm to be the extension of Frobenius norm of a matrix, i.e. $$
\|X\|^2=\sum_{i_1, \ldots i_d } X_{i_1, \ldots i_d}^2,
$$
and the tensor Chebyshev norm as $$
\|X\|_C=\max_{i_1, \ldots i_d } 
|X_{i_1, \ldots i_d}|.
$$

The $i$-th mode product of a tensor $X$ with a matrix $M \in \mathbb{R}^{m \times n_i}$ will give a new tensor $Y \in \mathbb{R}^{n_1\times \cdots n_{i-1}\times m \times n_{i+1} \ldots n_d}$ denoted by
$
Y=X \times_i M,
$
which is equivalent to $Y_{(i)}=M X_{(i)}$. With this notation a tensor $X$ with multi-linear rank $\mathbf{r}$ can then be represented by the Tucker decomposition 

\[
X=G \times_1 U_1 \times_2 U_2 \cdots \times_d U_d 
=G \times_{i=1}^d U_i.
\]
Here $G  \in \mathbb{R}^{r_1\times r_2 \times \cdots r_d} $ is the core tensor, and $U_i \in \mathbb{R}^{n_i \times r_i}$ are the factor matrices. In this paper, we let $U_i$ be orthonormal matrices, i.e. $U_i^T U_i=I_{r_i}.$ 

A quasi-optimal Tucker representation of a given tensor $X$ can be found using, e.g., the high-order singular value decomposition (HOSVD) \cite{tucker1966some,de2000multilinear}.  To this end, when discussing storage and computational complexity we, for simplicity, assume $n_i=n, r_i=r, i=1, \ldots d$. Under these assumptions the compressed Tucker format has storage requirement $O(r^d+dnr).$ When $r\ll n$ this results in significant storage savings compared to the full representation, which requires storage on the order of $O(n^d).$ 

\subsection{Review of matrix cross approximation}
In this subsection, we will review matrix cross approximation, which serves as a foundation for cross approximation of tensors. Low-rank approximation using matrix skeleton, pseudoskeleton, CUR factorization and cross approximation is a well-studied subject in numerical linear algebra \cite{goreinov1997theory,goreinov2001maximal,mahoney2009cur}. Given a matrix $G$, the cross approximation approximates $G$ by selected columns (indexed by $\mathcal{J}$) and rows (indexed by $\mathcal{I}$)
\[ 
G \approx G(:,\mathcal{J}) G(\mathcal{I},\mathcal{J})^+ G(\mathcal{I},:).
\]
The quality of this cross approximation critically hinges on the choice of the column and row indices, $\mathcal{I}, \mathcal{J}.$ Popular methods for choosing  $\mathcal{I}$ and $ \mathcal{J}$ include adaptive cross approximation \cite{bebendorf2000approximation}, randomized algorithms \cite{halko2011finding,chiu2013sublinear,cortinovis2025sublinear}, max volume \cite{goreinov2001maximal}, leverage score \cite{mahoney2009cur},
and DEIM (discrete empirical interpolation method) and its variant QDEIM (potentially with oversampling) \cite{sorensen2016deim,drmac2016new,peherstorfer2020stability}. The DEIM index selection methods apply DEIM on the leading left and right singular vectors of the matrix, and have been shown to offer robust performance in matrix CUR decomposition \cite{sorensen2016deim}. Recently, DEIM index selection has  been used for cross approximation in for low-rank methods for nonlinear PDEs and stochastic PDEs \cite{donello2023oblique,ghahremani2024deim,dektor2024collocation,ghahremani2024cross,naderi2024cur}. 

In \cite{appelo2025lraa}, we proposed an adaptive rank cross approximation based on DEIM index selection, called Cross-DEIM, to be used within lrAA. Given a user-specified tolerance Cross-DEIM outputs an approximate SVD of a matrix based on cross approximation. Cross-DEIM's index election is based on DEIM index selection \cite{sorensen2016deim} from approximate singular vectors and iteratively updates the approximate singular vector matrices, simultaneously enriching the column and row index sets. The iteration is stopped when consecutive iterates are close and a low-rank criteria is satisfied (this type of sampling strategy is called progressive sampling in \cite{xia2024making}). At the end of the iteration, a rank pruning process based on a user specified tolerance is carried out. We showed by numerical experiments that the method performed well for low-rank matrix approximation tasks   and particularly within lrAA with warm-start.  In the next subsection, we will propose   generalization of Cross-DEIM to Tucker tensors. 

\subsection{Cross$^2$-DEIM approximation for Tucker tensor}
\label{sec:cd}
  Cross approximations for Tucker tensors are generalizations of cross approximations for matrices with the objective to obtain low rank Tucker approximations by sampling part of the tensor. For a comprehensive overview of Tucker cross approximation algorithms, we refer to the review paper \cite{ahmadi2021cross} which surveys  various techniques including randomized and deterministic methods, and details their connections to cross approximations for matrices.

Our goal is to develop a fast Tucker cross approximation that has linear computational scaling with respect to $n$. Various methods in the literature satisfy this criteria, linear complexity in $n$, by performing fiber sampling. In \cite{oseledets2008tucker}, the authors proved the quasi-optimality of their Cross3D algorithm, which is based on a sequence of matrix cross approximation for the unfolding matrices. The computational cost of Cross3D is $O(n^2 r^2)$ for a  3D tensor. To achieve linear scaling in $n$, \cite{oseledets2008tucker} further proposed to conduct additional cross approximation for the slices which reduces the cost to $O(nr^4)$ for a 3D tensor. 
Two types of Fiber Sampling Tucker Decomposition (FSTD) were proposed in \cite{caiafa2010generalizing}. The first type, FSTD1, samples $O(r^{d-1})$ fibers in each mode, while FSTD2 is based on sampling $O(r)$ fibers in each mode. The paper \cite{caiafa2010generalizing} provides a practical algorithm for FSTD1 with cost $O(n r^4)$ for 3D tensors, but \cite{caiafa2010generalizing} does not provide an algorithms for the   faster FSTD2. Very recently, a Tucker cross approximation based on DEIM index selection was proposed in \cite{ghahremani2024deim}. The method resembles FSTD1, but uses a DEIM-based sampling; its computational cost is $O(nr^{2(d-1)}).$ Furthermore, the work \cite{ghahremani2024deim}  proposes an iterative procedure which allows adaptive rank approximation. Finally, the work in \cite{dolgov2021functional} propose a functional Tucker approximation for Chebyshev interpolation by iteratively refining the Chebyshev grids to achieve low cost in tensor entry sampling.

The goal of this work is to develop a practical algorithm for FSTD2.
For convenience of the readers, below we summarize the main theoretical results in \cite{oseledets2008tucker} and \cite{caiafa2010generalizing}, relating to FSTD2-type fiber sampling.

\begin{theorem}[Adapted from Cross3D \cite{oseledets2008tucker}]
\label{thm:ose}Suppose that a tensor $X$ admits a rank $(r_1, \ldots, r_d)$ Tucker approximation 
$$
X=G \times_1 U_1 \times_2 U_2 \cdots \times_d U_d +E,
$$
with $\|E\| \le \epsilon.$ Then there exists a selection of $r_i$ fibers in the $i$-th mode, which are represented as matrices $C_i \in \mathbb{R}^{n_i \times r_i}$ and a core tensor $G'\in \mathbb{R}^{r_1\times r_2 \times \cdots r_d},$ such that 
$$
X=G' \times_1 C_1 \times_2 C_2 \cdots \times_d C_d +E',
$$
where $\|E'\|_C \le f(r_1, \ldots r_d) \epsilon,$ and $f(\cdot)$ has polynomial dependence on its argument.
\end{theorem}
We remark that the proof of Theorem \ref{thm:ose} is constructive and is based on repeated use of the  matrix Cross2D algorithm.  

\begin{theorem}[Adapted from FSTD2  \cite{caiafa2010generalizing}]
\label{thm:fstd2}
Suppose that the tensor $X$ has multi-linear rank $(r, \ldots, r).$ Then, given a selection of $r$ fibers arranged as matrices $C_i \in \mathbb{R}^{n_i \times r}, i=1, \ldots d$ with rank $r,$ and given index sets $\mathcal{I}_i, i=1, \ldots d,$ such that the matrices $W_i=C_i(\mathcal{I}_i, :)$ are non-singular, we have
$$
X=G \times_1 C_1 \times_2 C_2 \cdots \times_d C_d,
$$
where
$$
G=W \times_1 (W_1)^{-1} \times_2 (W_2)^{-1} \cdots \times_d (W_d)^{-1},
$$
and $W=X(\mathcal{I}_1, \cdots, \mathcal{I}_d)$.
\end{theorem}

While not explicitly  stated  in the theorem, a key ingredient in FSTD2 is the use of two sets of indices. The first set is the one spelled out in Theorem \ref{thm:fstd2}: $\mathcal{I}_i \subset \mathcal{N}_i,$ which we call the main index set. The second set of indices is implied by the ``... given a selection of $r$ fibers ...''. This set is used to define the fibers to form the $C_i$ matrices in Theorems \ref{thm:ose} and \ref{thm:fstd2}. We denote those indices by $\mathcal{I}_{\neq i}$ and note that $\mathcal{I}_{\neq i} \subset \otimes_{j\neq i} \mathcal{N}_j$. We call $\mathcal{I}_{\neq i}$ the complement index set.

\subsubsection{Inner and adaptive Cross$^2$-DEIM algorithms}
Our algorithm Cross$^2$-DEIM (C2D for short), which we will now describe, provides a fast way to determine {\em both} index sets $\{\mathcal{I}_i, \mathcal{I}_{\neq i}, i=1 \ldots d \}.$ This method is generalized from Cross-DEIM and \cite{ghahremani2024deim} in the sense that DEIM is used for index selection of $\mathcal{I}_i$, but the key novelty lies in the selection of index sets 
$\mathcal{I}_{\neq i}.$ We present two versions of the methods, the inner algorithm, which contains the main ideas and the adaptive algorithm, which is an adaptive rank algorithm given a prescribed tolerance.

Starting with the inner version of C2D (denoted {\tt C2Di}), we first give an overview of the algorithm and then proceed to explain it in detail. This method requires an initial guess of the factor matrices $(U_i^{(0)},) i=1, \ldots d$. Then we use DEIM to select the first set of indices $\{\mathcal{I}_i, i=1 \ldots d\}$. That is, given an approximate singular vector matrix $U_i^{(0)}$ in each mode, DEIM$(U_i^{(0)})$ outputs an index set $\mathcal{I}_i$. The next step is to select the index sets $\{\mathcal{I}_{\ne i}, i=1 \ldots d\}.$ We note that the possible indices satisfy $\mathcal{I}_{\neq i} \subset \otimes_{j\neq i} \mathcal{N}_j$ but probing all of $\otimes_{j\neq i} \mathcal{N}_j$ with DEIM is too expensive. To achieve a faster method we make the \emph{assumption} that $\mathcal{I}_{\ne i} \subset \otimes_{j\neq i} \mathcal{I}_j,$ and select $\mathcal{I}_{\ne i}$ using a \emph{matrix cross approximation strategy}. 

In particular, for $W_{(i)}$, the $i$-th mode matricization of the core tensor $W=X(\mathcal{I}_1, \cdots, \mathcal{I}_d)$ we use DEIM on the leading $r_i$ right singular vectors of $W_{(i)},$ to  obtain the subset $\mathcal{I}_{\ne i}$ of  $r_i = |\mathcal{I}_{\ne i}|$ indices from $\otimes_{j\neq i} \mathcal{I}_j$. Note that since size of the core tensor is independent of $n$, the computational cost of finding $\mathcal{I}_{\ne i}$ is also independent of $n$. 

\begin{algorithm}[h!]
\begin{minipage}{0.9\linewidth}
\caption{{\tt [$G, U^{(1)}_1, \ldots, U^{(1)}_d $] = C2Di($X, U^{(0)}_1, \ldots U^{(0)}_d $)} \\Cross$^2$-DEIM approximation for Tucker tensor (inner algorithm) \label{alg:C2Di}}
 \begin{algorithmic}[1]
 \STATE {\bf Input:} Tensor  $X \in \mathbb{R}^{n_1\times n_2 \times \cdots n_d}$ and its approximate orthonormal factor matrices \mbox{$U_i^{(0)} \in \mathbb{R}^{n_i \times r_i}, i=1, \ldots d$.} 
\STATE {\bf Output:} Tucker decomposition $X \approx G \times_1 U^{(1)}_1 \times_2 U^{(1)}_2 \cdots \times_d U^{(1)}_d $ with $G \in \mathbb{R}^{r_1\times r_2 \times \cdots r_d},$ and orthonormal factor matrices $U_i^{(1)} \in \mathbb{R}^{n_i \times r_i}, i=1, \ldots d$. 
\\\hrulefill
\STATE $r_i = {\tt size}(U^{(0)}_{i},2), \, i = 1,\ldots,d$
\STATE $\mathcal{I}_{i} = {\tt QDEIM}(U^{(0)}_{i}), \, i=1, \ldots, d$ \COMMENT{QDEIM index selection for $\mathcal{I}_i$}
\STATE $\mathcal{I}^{\rm os}_{i} = {\tt gpode}(U^{(0)}_{i}, r_i+3), \quad i=1, \ldots, d$ \COMMENT{Oversampling index by 3}
\STATE $W=X(\mathcal{I}_1, \ldots, \mathcal{I}_d)$ \COMMENT{Cross tensor}
\STATE $W^{\rm os}=X(\mathcal{I}^{\rm os}_1, \ldots, \mathcal{I}^{\rm os}_d)$ \COMMENT{Oversampled cross tensor}
\FOR{$i = 1, \ldots, d$}
\STATE $W_{(i)}=\verb+core_reshape+(W, r_i, R_{\ne i})$ 
\STATE $[Z_i,\sim,\sim]={\tt svd}(W_{(i)}^T,{\verb+"econ"+})$ \COMMENT{Compute $r_i$ sing. vectors for second index set}
\STATE  $\mathcal{J}_{\neq i} = {\tt QDEIM}(Z_{i})$ \COMMENT{Index selection for $\mathcal{I}_{\neq i}$}
\STATE $\mathcal{I}_{\neq i} = J(\mathcal{J}_{\neq i})$ \COMMENT{Convert indices  for $\mathcal{J}_{\neq i}$ to $\mathcal{I}_{\neq i}$}
\STATE $\hat{W}_{(i)}=\verb+reshape+(X([\mathcal{I}_{\neq i}]_1, \ldots, \mathcal{N}_i ,\ldots, [\mathcal{I}_{\neq i}]_d),n_i,r_i)$ \COMMENT{Sample fibers }
\STATE $[U_{i}^{(1)},\sim,\sim]={\tt svd}(\hat{W}_{(i)},{\verb+"econ"+})$ \COMMENT{matrices and orthogonalize}
\ENDFOR
\STATE $G = W^{\rm os} \times_1 [U_{1}^{(1)}(\mathcal{I}^{\rm os}_1,:)]^+ \times_2  [U_{1}^{(1)}(\mathcal{I}^{\rm os}_2,:)]^+ \cdots  \times_d [U_{d}^{(1)}(\mathcal{I}^{\rm os}_d,:)]^+$

\STATE Return $G,  U^{(1)}_1, \ldots U^{(1)}_d$.
\end{algorithmic}
\end{minipage}
\end{algorithm}

The procedure {\tt C2Di}, outlined above, is detailed in Algorithm \ref{alg:C2Di}, which we now explain line by line. Lines 4 - 8 are self explanatory (the algorithms \verb+QDEIM+ and \verb+gpode+ can be found in \cite{drmac2016new} and \cite{peherstorfer2020stability}). To explain line 10 it is convenient to introduce the linear ordering of the $\prod_{i=1}^d r_i$ elements in $W = X(\mathcal{I}_1, \cdots, \mathcal{I}_d)$ as 
\begin{equation} \label{eq:linorder}
J(k_1,k_2,\ldots,k_d) = k_1 + \sum_{l=2}^{d}  \prod_{s=2}^{l} (k_s-1)r_{s-1}.
\end{equation} 
We also denote by $R_{\ne i}$ the $d$ integers $R_{\neq i} = \prod_{j = 1, j \ne i}^d r_j,\, i = 1\ldots,d$.
Then by the dimension $i$ dependent function \verb+core_reshape+$(W, r_i, R_{\ne i})$ used in line 10 of Algorithm \ref{alg:C2Di} we mean the reordering of the elements in $W$ into a matrix of size $r_i \times R_{\ne i}$ with rows $k_i =  1,2 \ldots, r_i$ and columns $j = 1,2,\ldots, R_{\ne i}$, i.e.
\[
W_{(i)}(k_i,j) = X(\mathcal{I}_1(k_1),\mathcal{I}_2(k_2),\ldots,\mathcal{I}_d(k_d)).
\]  
Here the arguments in the right hand side $(k_1,k_2,\ldots,k_d)$ follow the linear ordering (\ref{eq:linorder}). After computing the $R_{\ne i} \times r_i $ leading the left singular vector matrix $Z_{i}$ of of $W_{(i)}^T$ (line 11) we can use DEIM on $Z_{i}$ get an index selection $\mathcal{J}_{\ne i} = {\tt QDEIM}(Z_{i})$ from the index set $\{1,2,\ldots, R_{\ne i}\}$ (line 12). Using the linear ordering (\ref{eq:linorder}), this index set can then be mapped back to $r_i$ $(d-1)$-tuples of indices from $\{\mathcal{I}_1,\ldots,\mathcal{I}_{i-1},\mathcal{I}_{i+1},\ldots, \mathcal{I}_d\}$. This selection of $r_i$ tuples form our index set $\mathcal{I}_{\ne i}$. On line 13 in Algorithm \ref{alg:C2Di}  we use the shorthand $\mathcal{I}_{\ne i} = J(\mathcal{J}_{\ne i})$ to denote the mapping from $\mathcal{J}_{\ne i}$ to $\mathcal{I}_{\ne i}$. Next, in line 14, we use the index set $\mathcal{I}_{\ne i}$ to form the $n_i \times r_i$ matrix $\hat{W}_{(i)}$ with elements 
\begin{multline*}
\hat{W}_{(i)}(j,k) = X([\mathcal{I}_{\ne i}(k)]_1,\ldots,[\mathcal{I}_{\ne i}(k)]_{i-1},j,[\mathcal{I}_{\ne i}(k)]_{i+1},\ldots,[\mathcal{I}_{\ne i}(k)]_d), \\ j = 1,\ldots,n_i, \ \ k = 1,\ldots, r_i.
\end{multline*} 
Here we use the notation $[\mathcal{I}_{\ne i}(k)]_l$ to denote the $k$th element in the index set in dimension $l$. These matrices are then orthogonalized (line 15) into the factor matrices, $U_{i}^{(1)}$,  returned as outputs. Finally, in line 17, we form the output core matrix using the factors $U_{i}^{(1)}$ and the oversampled cross matrix $W^{\rm os}$ \[
G = W^{\rm os} \times_1 [U_{1}^{(1)}(\mathcal{I}^{\rm os}_1,:)]^+ \times_2  [U_{1}^{(1)}(\mathcal{I}^{\rm os}_2,:)]^+ \cdots  \times_d [U_{d}^{(1)}(\mathcal{I}^{\rm os}_d,:)]^+.
\]
Here $A^{+}$ denotes the formal application of the pseudo-inverse. In practice, the oversampled problem is solved using least squares. We note that the output approximate tensor has multi-linear rank $(r_1, \ldots, r_d).$

\begin{algorithm}[h!]
\begin{minipage}{0.9\linewidth}
\caption{{\tt [$G, U_1, \ldots, U_d $] = C2D($X, U^{(0)}_1, \ldots U^{(0)}_d $,{\rm TOL},$q$,${\rm iter}_{\rm max}$)} \\Cross$^2$-DEIM approximation for Tucker tensor (adaptive algorithm) \label{alg:C2Da}}
 \begin{algorithmic}[1]
 \STATE {\bf Input:} Tensor  $X \in \mathbb{R}^{n_1\times n_2 \times \cdots n_d}$ and its approximate orthonormal factor matrices \mbox{$U_i^{(0)} \in \mathbb{R}^{n_i \times r_i}, i=1, \ldots d$, tolerance {\rm TOL}, rank increase $q$, max iterations ${\rm iter}_{\rm max}$.} 
\STATE {\bf Output:} Tucker decomposition $X \approx G \times_1 U_1 \times_2 U_2 \cdots \times_d U_d $ with $G \in \mathbb{R}^{r_1\times r_2 \times \cdots r_d},$ and orthonormal factor matrices $U_i \in \mathbb{R}^{n_i \times r_i}, i=1, \ldots d$. 
\\\hrulefill

\STATE $ \mathcal{I}_{i}^{0} = \emptyset, i = 1,\ldots,d $ \COMMENT{Start from empty index sets}
\FOR{$k = 1,2,\ldots, {\rm iter}_{\rm max}$}
\FOR{$i = 1, \ldots, d$}
\STATE $r_i = {\tt size}(U^{(k-1)}_{i},2)$
\STATE $\mathcal{I}_{i}^{\ast} = {\tt QDEIM}(U^{(k-1)}_{i})$
\STATE $\mathcal{I}_{i} = \mathcal{I}_{i}^k$ = \verb+increase+$(\mathcal{I}_{i}^{\ast},\mathcal{I}_{i}^{k-1},\mathcal{N}_i,q)$ \COMMENT{Grow {\tt QDEIM} index set by $q$ elements}
\STATE $\mathcal{I}^{\rm os}_{i} = {\tt gpode}(U^{(k-1)}_{i}, r_i+3)$ 
\ENDFOR
\STATE $W=X(\mathcal{I}_1, \ldots, \mathcal{I}_d)$,  $W^{\rm os}=X(\mathcal{I}^{\rm os}_1, \ldots, \mathcal{I}^{\rm os}_d)$
\FOR{$i = 1, \ldots, d$}
\STATE $W_{(i)}=\verb+core_reshape+(W, r_i, R_{\ne i})$ 
\STATE $[Z_i,\sim,\sim]={\tt svd}(W_{(i)}^T,{\verb+"econ"+})$ 
\STATE  $\mathcal{J}_{\neq i} = {\tt QDEIM}(Z_{i})$
\STATE $\mathcal{I}_{\neq i} = J(\mathcal{J}_{\neq i})$ 
\STATE $\hat{W}_{(i)}=\verb+reshape+(X([\mathcal{I}_{\neq i}]_1, \ldots, \mathcal{N}_i ,\ldots, [\mathcal{I}_{\neq i}]_d),n_i,r_i)$
\STATE $[U_{i}^{(1)},\sim,\sim]={\tt svd}(\hat{W}_{(i)},{\verb+"econ"+})$ 
\ENDFOR
\STATE $G^{(k)} = W^{\rm os} \times_1 [U_{1}^{(k)}(\mathcal{I}^{\rm os}_1,:)]^+ \times_2  [U_{1}^{(k)}(\mathcal{I}^{\rm os}_2,:)]^+ \cdots  \times_d [U_{d}^{(k)}(\mathcal{I}^{\rm os}_d,:)]^+$
\STATE $\sigma_{{\rm min},i} = \verb+min_sing_val+(G^{k}_{(i)}), \, i = 1,\ldots,d$ \COMMENT{Compute min sing. val.} 
\IF{$\max \left(\| G^{k} \times_{i=1}^d U_i^{k}  - G^{k-1} \times_{i=1}^d U_i^{k-1} \|, \max_{j=1}^d \left( \sigma_{{\rm min},j}\right) \right) < {\rm TOL}$  }
\STATE $[G,U_1, \ldots U_d] = \verb+HOSVD+(G^{k},U_1^{k}, \ldots, U_d^{k},{\rm TOL})$
\STATE Exit iteration 
\ENDIF
\ENDFOR
\STATE Return $G,  U_1, \ldots, U_d$.
\end{algorithmic}
\end{minipage}
\end{algorithm}

We now describe the adaptive version of Cross$^2$-DEIM (C2D) described in Algorithm \ref{alg:C2Da}, which gives a cross approximation with adaptive rank up to prescribed tolerance {\rm TOL}. The idea is to improve the guess of the factor matrices $U_i^{(0)}$ by an iterative procedure. The core elements of the C2D are the same as for C2Di and only line 9 and 22-25 need further explanation. Starting with line 9, the procedure \verb+increase+ takes the most recent index set $\mathcal{I}^{\ast}_i$ and adds $q \ge 0$ more {\bf unique} elements from $\mathcal{I}^{k-1}_i$  (if available). If $q>0$ the rank increases monotonically by $q$ in each iteration. Choosing $q$ is a balance between memory limit and iteration count. A small value of $q$ will ensure moderate growth of intermediate ranks,  however, it will also mean more iteration steps are needed for higher target rank.  In all numerical experiments in this paper, we fix $q=4$ and find that to be a good choice. We note that in Cross-DEIM for matrix \cite{appelo2025lraa}, we did not have this bound on the rank increase, instead we always merge the old and new index set. Due to the high dimensionality of tensors, we find it inefficient to do so because of the rapid rank growth.
We also note that if we set $q=0,$ this will turn the method into a fixed rank algorithm.

On line 22 we compute the smallest singular values of the $d$ unfoldings of the $k$th iterate core tensor. This is to check if the low rank condition is satisfied. On line 23 we use these together with the norm of the difference between the two last iterates to break out of the iteration. We break when all the smallest singular values and  the difference between the two last iterates, in norm, are smaller than the user provided tolerance. Finally, on line 24 we perform a final truncation using the high order SVD with absolute tolerance {\rm TOL}.

\subsubsection{Discussion}
The computational cost of Algorithm \ref{alg:C2Di} is $O(dn r^2+r^{d+1})$ with $O(dnr+r^d)$ evaluations of tensor entries, which is  lower than the existing ones in the literature  \cite{ahmadi2021cross}. The computational cost of Algorithm \ref{alg:C2Da} depends on the iteration number, which stays $O(1)$ in practice, and the  intermediate ranks. In particular, C2D excels when a  warm-start strategy can be incorporated. This will reduce the iteration numbers significantly. Namely, for parametric tensor approximations, we can input the initial guess of factor matrices $U^{(0)}_i$ from a close and already known approximation to reduce iteration number. This is very useful in many computational tasks including Tucker-AA to be described in the next section.

\section{Anderson Acceleration in Tucker tensor format}
\label{sec:tuckeraa}
The AA method, \cite{anderson1965iterative, saad2024acceleration}, seeks to accelerate a Picard fixed-point iteration
\begin{equation} \label{eq:Picard}
X_{k+1} = H(X_k).
\end{equation}
It does so by using a weighted sum of the previous $\min(k,\hat{m})+1$ iterates and residuals (here $\hat{m}$ is called the window size or memory parameter). The properties of the full rank version of AA has been studied in terms of convergence and efficient formulations \cite{walker2011anderson,toth2015convergence,evans2018proof,zhang2018globally,pollock2019anderson,de2022linear,rebholz2023effect,de2024anderson} and it is also known that it has close connections to GMRES and multisecant methods \cite{walker2011anderson,fang2009two,lin2013elliptic}. 

In \cite{appelo2025lraa} we adopted AA to low-rank matrix format by replacing all linear and nonlinear operations by equivalent low-rank linear and nonlinear operations. Specifically, the unknown and intermediate iterates were stored as SVDs and the linear and nonlinear operations were carried out with   Cross-DEIM, which can be thought of as a matrix specialization of Cross$^2$-DEIM.   
The use of a cross approximation within the AA framework is naturally equipped with two rank limiting operations. First, the finite window size of AA is naturally limiting the rank growth of intermediate iterates, second the truncation in a cross approximation also controls the rank. A third mechanism to further control intermediate rank inflation is to use some type of scheduling for the tolerance used in cross approximation. Following \cite{appelo2025lraa} we simply set the tolerance to be proportional to the current residual.  We demonstrated, computationally, this often led to a monotonically increasing behavior of the intermediate solution rank throughout the iteration. In this section, we will present  Tucker-AA, the extension of lrAA to Tucker tensor format, as a solver for nonlinear tensor equations.

\begin{algorithm}[ht!]
\begin{minipage}{0.9\linewidth}
\caption{Tucker-AA for nonlinear tensor equation $H(X)=X.$  \label{alg:TAA} }
  \begin{algorithmic}[1]
  \STATE {\bf Input:} Initial guess $X_0 = G^{(0)} \times_{i=1}^d U^{(0)}_i$, window size $\hat{m} \ge 1,$ scheduling parameter $\theta$, tolerance ${\rm TOL}$, C2D max iter  ${\rm iter}_{\rm max}$, C2D increase $q$, max rank $r_{\rm max}$, truncation parameter $\epsilon_F$, initial C2D truncation level $\epsilon_{H}^0$.  
     \STATE {\bf Output:} Approximate solution $X_k = G^{(k)} \times_{i=1}^d U^{(k)}_i$ to the fixed point problem $H(X)=X$ in its Tucker tensor form.
     \\\hrulefill
  \STATE $\epsilon_{H} = \epsilon_{H}^0$. \COMMENT{Choose $\epsilon_{H}^0$ so that $H_0$ has low rank.}

\STATE $X_1 = H_0=\textsf{C2D}(H(X_0),U^{0}_1,\ldots,U^{0}_d,\epsilon_{H},q,{\rm iter}_{\rm max})$. 
\STATE   $\rho_0 = \| H_0 - X_0 \|$. \COMMENT{Compute residual}
\FOR{$k = 1, 2, \ldots$  }
\STATE 
$
H_k=\textsf{C2D}(H(X_k),U^{k}_1,\ldots,U^{k}_d,\epsilon_{H},q,{\rm iter}_{\rm max}).$
\STATE   $\rho_k = \| H_k -X_k\|$. \COMMENT{Compute residual}
\STATE  $\hat{m}_k = \min(\hat{m},k)$.
\STATE Set $D_k = (\Delta F_{k-\hat{m}_k},\ldots,\Delta F_{k-1})$, where 
\STATE $\Delta F_i =  \mathcal{T}^{\rm round}_{\epsilon_F,r_{\rm max}}(F_{i+1}-F_i$) and \mbox{$F_i =  \mathcal{T}^{\rm round}_{\epsilon_F,r_{\rm max}}(H_i - X_i)$}.
\STATE Solve the least square problem
\begin{equation*}
\gamma^{(k)} = {\rm arg}\!\!\min_{v\in\mathbb{R}^{\hat{m}_k}}\|F_k - \sum_{j=0}^{\hat{m}_k-1}  [D_k]_j v_j\|, \qquad \gamma^{(k)} = (\gamma^{(k)}_0,\ldots,\gamma^{(k)}_{\hat{m}_k-1}).
\end{equation*}
\STATE  Compute the new iterate \COMMENT{This sum can be done with C2D if $\hat{m}_k$ is large.}
\small
\[ 
X_{k+1}  = \mathcal{T}^{\rm round}_{\epsilon_{H},r_{\rm max}}(H_k -  \sum_{i=0}^{\hat{m}_k-1}  \gamma_i^{(k)} \left[H_{k-\hat{m}_k+i+1} - H_{k-\hat{m}_k+i} \right] ).
\]
\normalsize 
\STATE Set $\epsilon_H =  \theta \rho_k$. \COMMENT{Update the truncation tol.} 
\IF {$\rho_k < {\rm TOL}$} 
\STATE Exit and return $X_{k+1}$.
\ENDIF
\ENDFOR
    \end{algorithmic}
\end{minipage}
\end{algorithm}
\subsection{Tucker-AA}
 The Tucker-AA algorithm is described in Algorithm \ref{alg:TAA}.  The algorithm outputs an approximate solution to the fixed point problem \eqref{eq:nle} in its  Tucker tensor form $X_k = G^{(k)} \times_{i=1}^d U^{(k)}_i$, with factors $U_i^{(k)}\in \mathbb{R}^{n_i \times r_i}$, and core $G^{(k)} \in \mathbb{R}^{r_1\times \cdots \times r_d}$, with $r\le r_{\rm max}$. The output Tucker tensor meets the  tolerance in the sense that $\|H(X_k)-X_k\| \le {\rm TOL}$ as specified by the input parameter ${\rm TOL}$. We now describe the critical parts of the algorithm in some detail. 

\subsubsection{Linear and nonlinear rank truncations}
First we note that throughout the algorithm quantities like $X_k$ and $H_k$ should be interpreted as their Tucker tensor representation, we never form the full tensors, but always operate on their factors. The core component of the algorithm is the C2D method as defined in Algorithm \ref{alg:C2Da}. When operating on a nonlinear function (as in lines 5 and 8) we must use C2D but note that when operating on a linear combination of Tucker tensors (as in line 14) it is also possible to use a standard truncated sum of tensors (see for example Exercise 5.3 in combination with Algorithm 6.3 in \cite{koldabook}). The former will typically be more efficient when the window size is large. We denote this truncation operation by $\mathcal{T}^{\rm round}_{\epsilon,r_{\rm max}}(\cdot)$ and use it in line 12 to compute differences of residuals as needed to form the linear least squares problem (lines 11-13). We note that the tolerance $\epsilon_{F}$ should be chosen to be very small in order to avoid ill-conditioning. We recommend to set it to $10^{-12}$ or smaller. Note that in all applications of C2D, we use   the factor matrices from the previous iterate to warm-start the cross approximation. 

\subsubsection{Scheduling of tolerance}
The parameter $\epsilon_H$ is a tolerance to control the error in the C2D approximation. It serves a similar role as the truncation tolerance in the (truncated) HOSVD. 

The rationale for scheduling, the process of adjusting $\epsilon_H$ throughout the iteration, is to make the iterates gradually (and preferably monotonically) increase in rank until the AA iteration terminates. As mentioned above, we use a simple strategy as in \cite{appelo2025lraa} by setting
\begin{equation}
\label{eq:simps}
\epsilon_H   = \rho_k \theta.
\end{equation} 
Here the quantity $\rho_k$ is a residual used to determine when to stop the AA iteration and  $\theta$ is a user specified scheduling parameter. 

An alternative scheduling was proposed in \cite{bachmayr2017iterative} for fixed point iterations stemming from Richardson iteration for elliptic problem by hierarchical Tucker tensor, which has theoretical footings when used with soft thresholding. A third possible strategy is to set $\epsilon_H = 0$ and turn Algorithm \ref{alg:TAA} into a fixed rank method by using the $r_{\rm max}$ option. For brevity, in this paper, we focus on rank-adaptive version and do not experiment with fixed-rank Tucker-AA. In our numerical experiment, we always disable $r_{\rm max}$ by setting it to be the $\min_{i = 1}^d n_i.$

\subsubsection{Low-rank solution to the least squares problem}
Finally, line 13 in Algorithm \ref{alg:TAA} finds the solution to the least squares minimization problem using Algorithm \ref{alg:leastsquareT} to update the solution as a linear combination of previous $H_k$. Note that the Algorithm \ref{alg:leastsquareT} only uses low rank operations. For a memory parameter $\hat{m}$ and a problem with $n_i = n$ and factor ranks \mbox{$r_i^{j}= {\rm rank}(U_i^j)) = r$} the cost of the algorithm is $\mathcal{O}(dn (r\hat{m})^2 + (r\hat{m})^3)$.

 \begin{algorithm}[h]
 \begin{minipage}{0.9\linewidth}
 \caption{Computing the least squares solution minimizing \\
 $\| \sum_{j=1}^{s} \gamma_j G^j \times_{i=1}^d U^j_i - G^B \times_{i=1}^d U^B_i \|$	\label{alg:leastsquareT}}
    \begin{algorithmic}[1]
      \STATE {\bf Input:} Tucker tensors $G^j \times_{i=1}^d U^j_i, j=1, \ldots s$,   $G^B \times_{i=1}^d U^B_i$. The ranks of the factors $U^{j}_i$ are denoted $r_{i}^j$.
      \STATE {\bf Output:} $\gamma_j, j=1, \ldots s$
      \\\hrulefill
     \FOR{$i = 1, \ldots, d$}
      \STATE Form $U_i=[U^1_i,\dots ,U^s_i]$
      \STATE $[Q_i,R_i,\Pi_i] = \verb+qr+(U_i,\verb+"econ"+)$ \COMMENT{Perform a column pivoted QR} 
      \STATE $c_i = [0, r_i^1, r_i^1+r_i^2, \ldots, \sum_{j=1}^s r_i^j]$    \COMMENT{$c_i$ is a vector of the cumulative sum}
\STATE      \COMMENT{of the ranks in dimension $i$}
    \ENDFOR    
    \STATE $b = {\sf vec}(G^B \times_{i=1}^d (Q_i^T U^B_i))$
    \FOR{$j = 1, \ldots, s$}
    \FOR{$i = 1, \ldots, d$}
    \STATE $Z = R_i \Pi_i^T$
    \STATE $W_i = Z(:,c_i(j)+1:c_i(j+1))$
    \ENDFOR
    \STATE	 $A(:,j) = {\sf vec}(G^j \times_{i=1}^d (W_i))$ \COMMENT{Form columns of the LS problem}
    \ENDFOR
    \STATE $\gamma = A^{+}b$ \COMMENT{Solve the least squares system (using QR or SVD)}
    \STATE Return $\gamma$
\end{algorithmic} 
\end{minipage}
\end{algorithm} 

\section{Numerical examples}
\label{sec:num}
In this section, we present numerical examples that demonstrate the methods: Cross$^2$-DEIM and Tucker-AA. 

\subsection{Cross$^2$-DEIM approximation for problems with cold-start}
In this experiment we consider the approximation of $n_1 \times n_2 \times n_3$ tensors whose elements are given by one of the four functions 
\begin{align*}
& X_1(i_1,i_2,i_3) = \frac{1}{i_1+i_2+i_3},  && X_{2}(i_1,i_2,i_3) = e^{-\left(x_{1}(i_1) x_{2}(i_2) x_{3}(i_3) \right)^2}, \\
&X_3(i_1,i_2,i_3) = (i_1^3+i_2^3+i_3^3)^{-\frac{1}{3}}, && X_4(i_1,i_2,i_3) = (i_1^5+i_2^5+i_3^5)^{-\frac{1}{5}}.  
\end{align*}
For $X_1$ and $X_2$ we take $n_1 = n_2 = n_3 =100$ and for $X_3$ and $X_4$ we take $n_1 = n_3 = 300$ and $n_2 = 400$. In $X_2$ we have $x_k(i_l) = -1 + 2\frac{l-1}{n_k-1}$.  By cold-start, we mean that each factor matrix used for the initial guess is a random vector with norm one.

\graphicspath{{figures}}
\begin{figure}[]
\begin{center}
\includegraphics[width=0.48\textwidth,trim={0.0cm 0.0cm 0.0cm 0.0cm},clip]{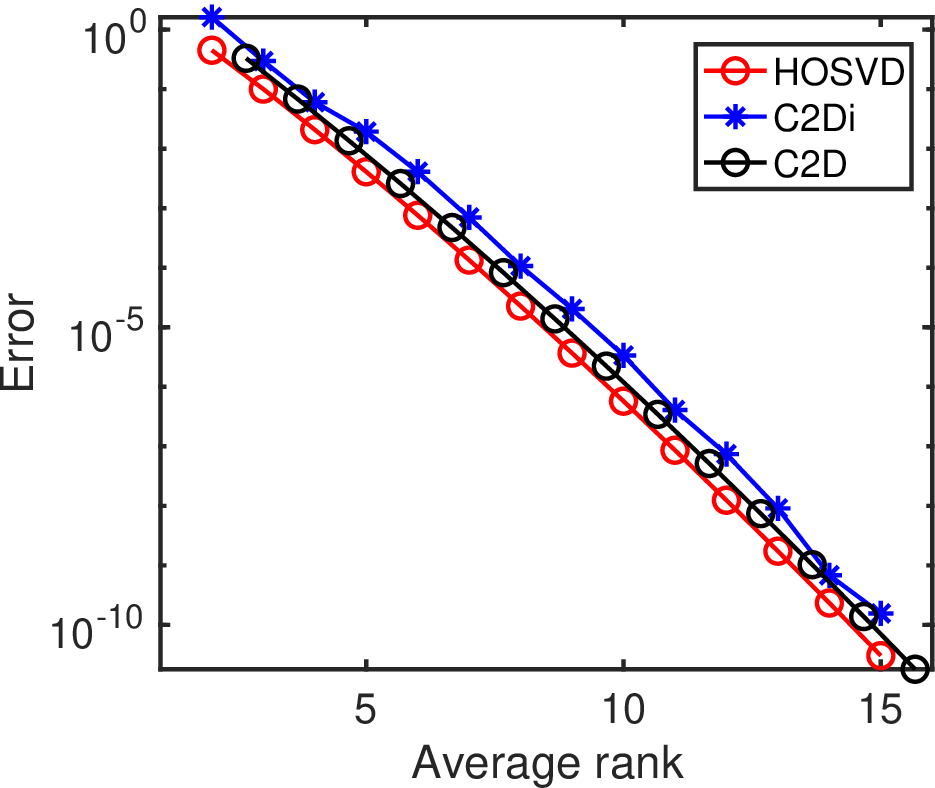}
\includegraphics[width=0.48\textwidth,trim={0.0cm 0.0cm 0.0cm 0.0cm},clip]{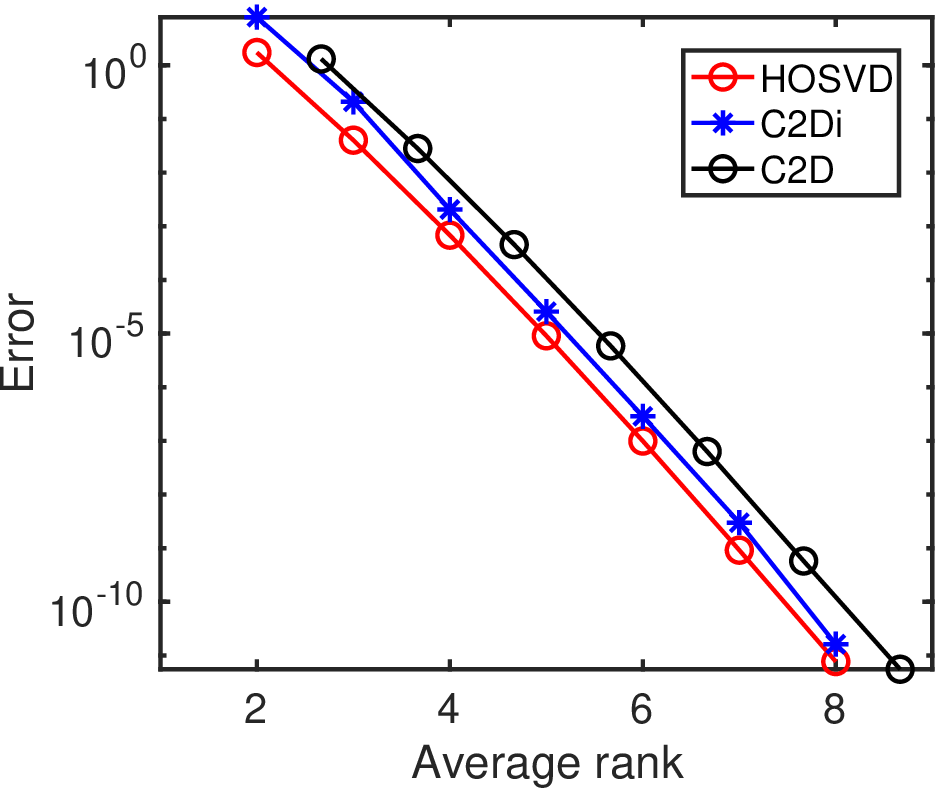}
\includegraphics[width=0.48\textwidth,trim={0.0cm 0.0cm 0.0cm 0.0cm},clip]{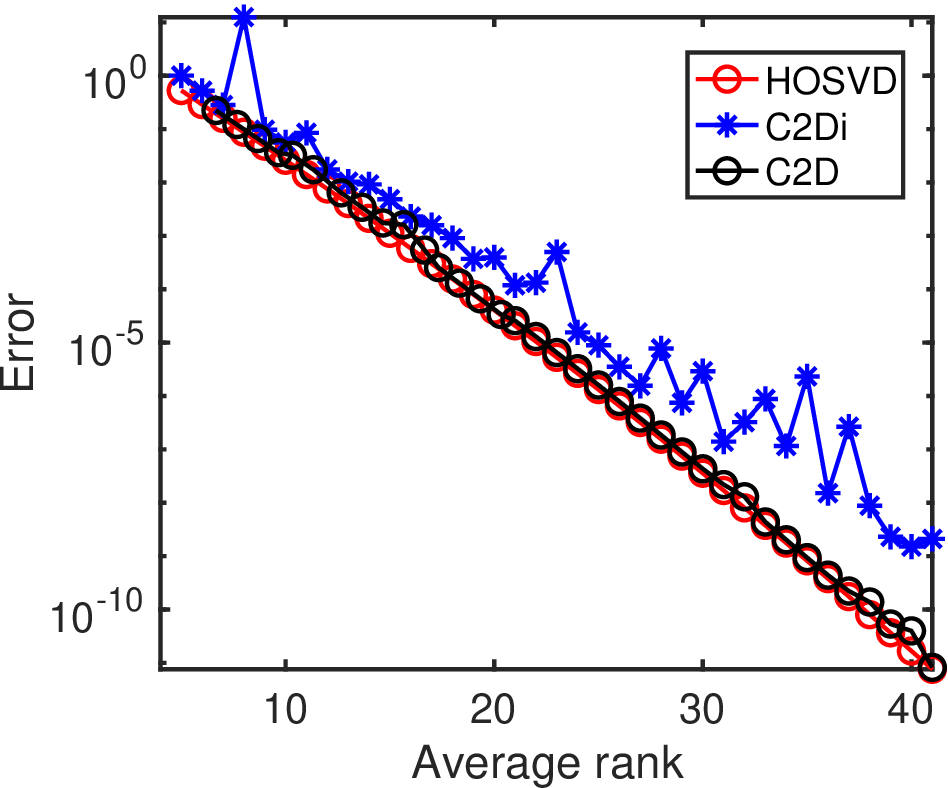}
\includegraphics[width=0.48\textwidth,trim={0.0cm 0.0cm 0.0cm 0.0cm},clip]{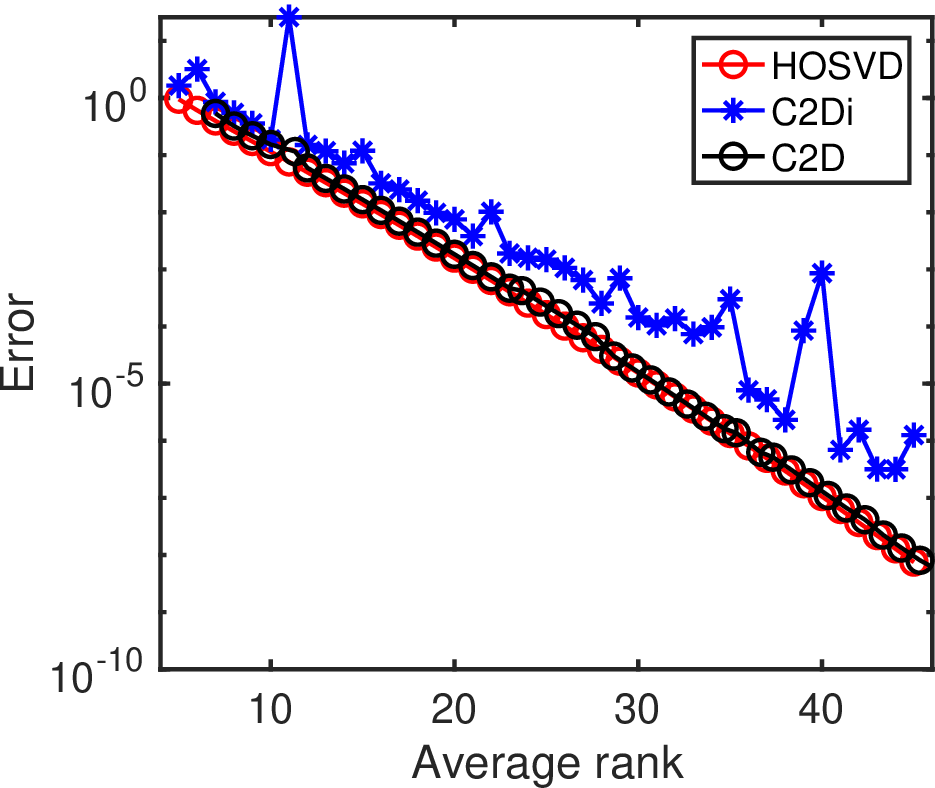}
\caption{Errors as a function of the average rank of the approximation of tensors $X_1, X_2, X_3, X_4$ (top left to bottom right) for HOSVD, C2Di and C2D. See the text for details. \label{fig:functions_error}}
\end{center}
\end{figure}

In Figure \ref{fig:functions_error}, in red, we display the absolute errors (in Frobenius norm) for the Tucker approximation computed using the HOSVD with a fixed rank in all three dimensions. We also display, in blue, the errors in the approximation obtained with C2Di using the HOSVD approximation as input. Finally, in black, we display the errors obtained using C2D with the tolerance selected to the HOSVD errors at varying ranks and with a rank 1 random tensor as initial guess. Since C2D may return a tensor with different ranks in the three dimensions we display the error as a function of the average rank (the average taken over the three dimensions). The results in Figure \ref{fig:functions_error} indicate that the performance of C2Di and C2D for approximations of $X_1, X_2$. However, for $X_3, X_4,$ we can see that C2Di does not offer robust decrease of error with respect to rank growth. This is due to the complement index selection  where only the sub-sampled tensors are looked at. However, with the iterative procedure in C2D, we are able to achieve  robust performance for the error behavior.    

\begin{figure}[]
\begin{center}
\includegraphics[width=0.48\textwidth,trim={0.0cm 0.0cm 0.0cm 0.0cm},clip]{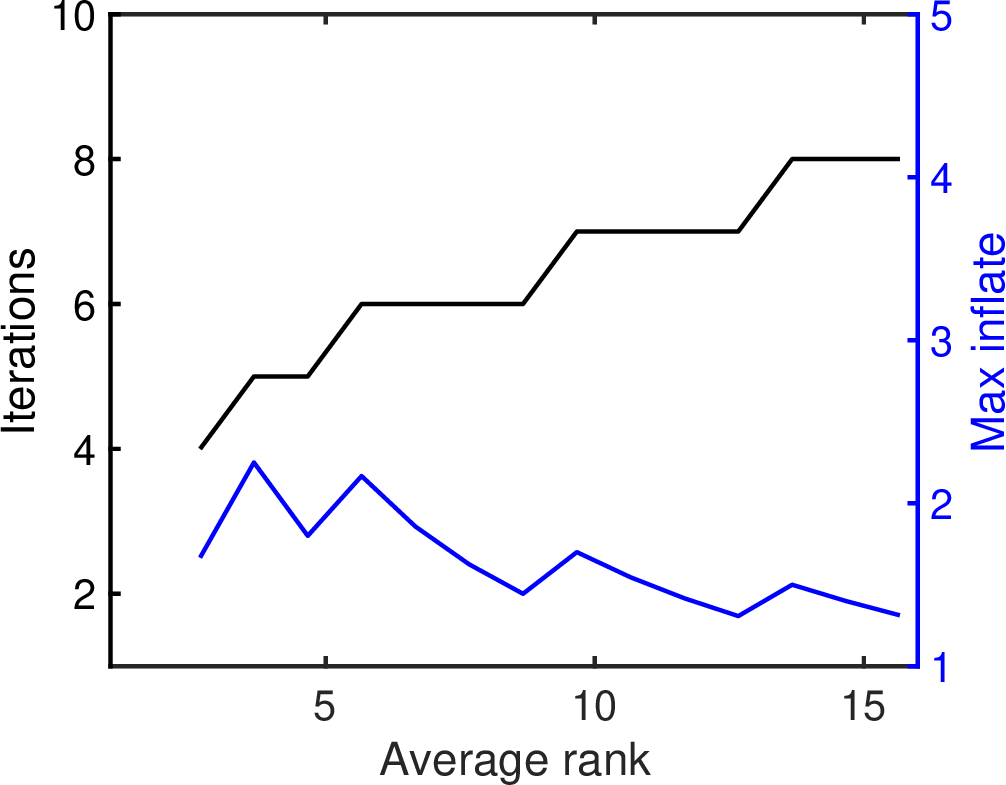}
\includegraphics[width=0.48\textwidth,trim={0.0cm 0.0cm 0.0cm 0.0cm},clip]{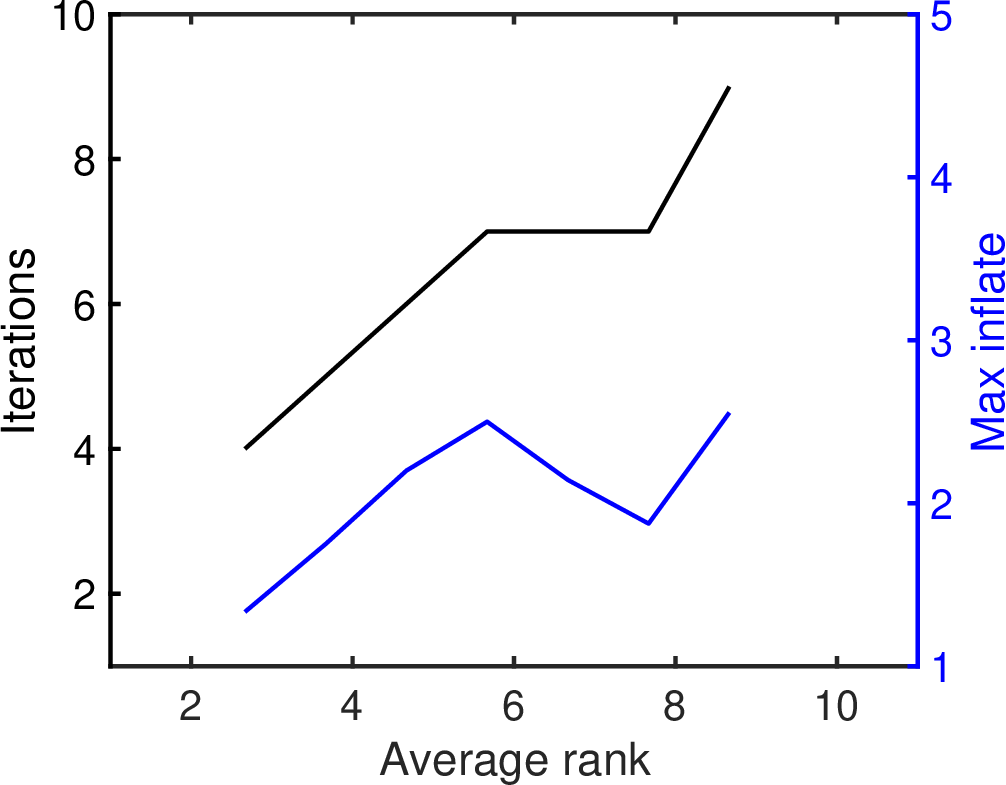}
\includegraphics[width=0.48\textwidth,trim={0.0cm 0.0cm 0.0cm 0.0cm},clip]{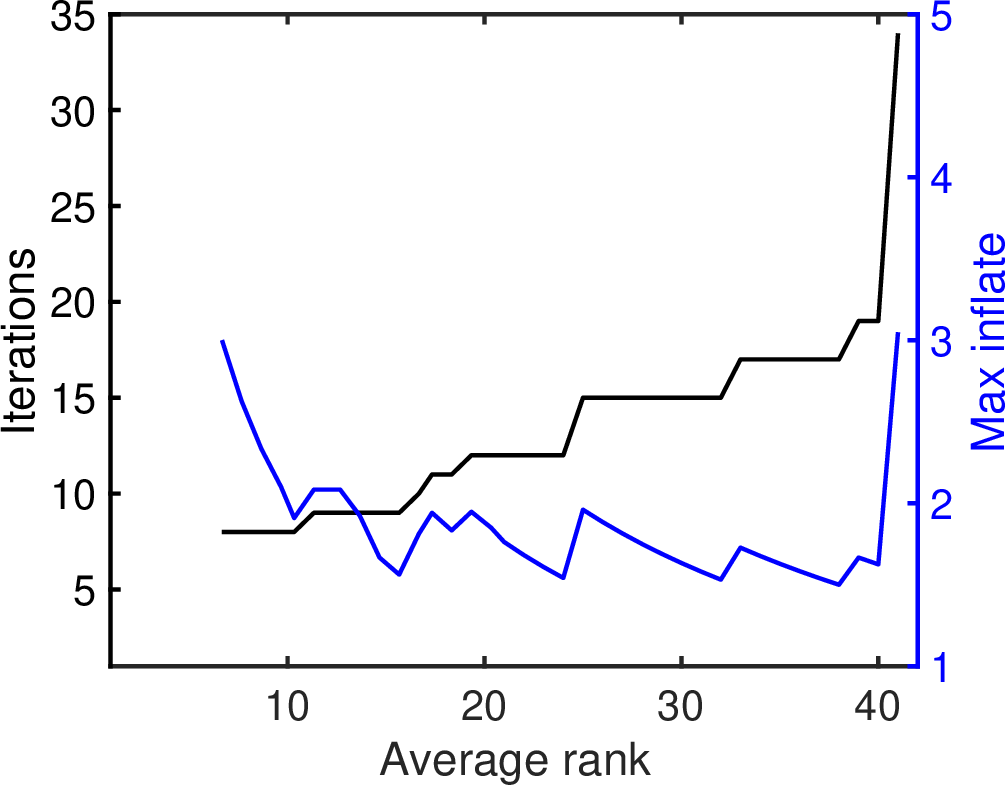}
\includegraphics[width=0.48\textwidth,trim={0.0cm 0.0cm 0.0cm 0.0cm},clip]{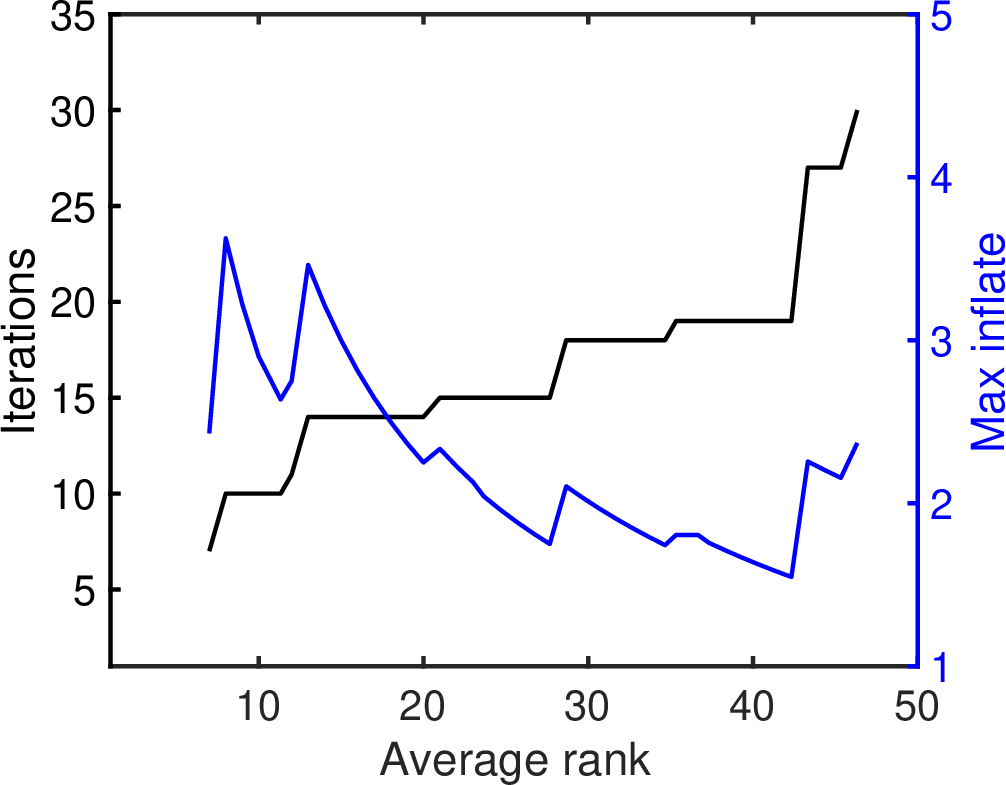}
\caption{Iterations and inflation as a function of the average rank of the approximation of tensors $X_1, X_2, X_3, X_4$ (top left to bottom right) for C2D. See the text for details. \label{fig:functions_rank}}
\end{center}
\end{figure}

In Figure \ref{fig:functions_rank}, in black, we display the number of iterations used by C2D for the four tensors and as a function of the average rank. We also display maximum inflation, which we define as the ratio of the maximum average rank during the C2D iteration and the average rank of the converged solution. The figures reveal that C2D converges in relatively few iterations even when started from a random initial guess and that for most tolerances the maximum inflation is around two or less. For the approximations of $X_3$ and $X_4,$ we do notice an increase in iteration count for tight tolerances. We remark that the parameter $q$, which is the rank increase in each iteration (here chosen as $q=4$), prevents the number of iterations to be lower than ``target rank'' divided by $q$.  

\subsection{Cross$^2$-DEIM approximation for problems with warm-start}
\begin{figure}[]
\begin{center}
\includegraphics[width=0.45\textwidth,trim={0.0cm 0.0cm 0.0cm 0.0cm},clip]{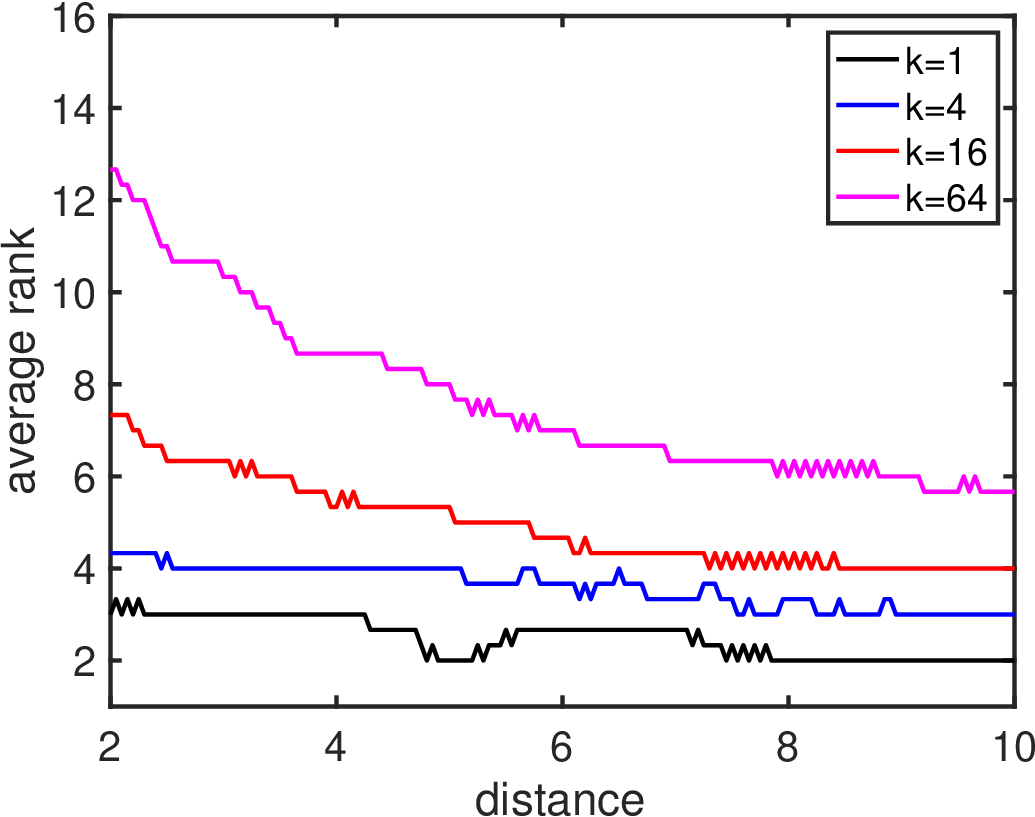} 
\includegraphics[width=0.45\textwidth,trim={0.0cm 0.0cm 0.0cm 0.0cm},clip]{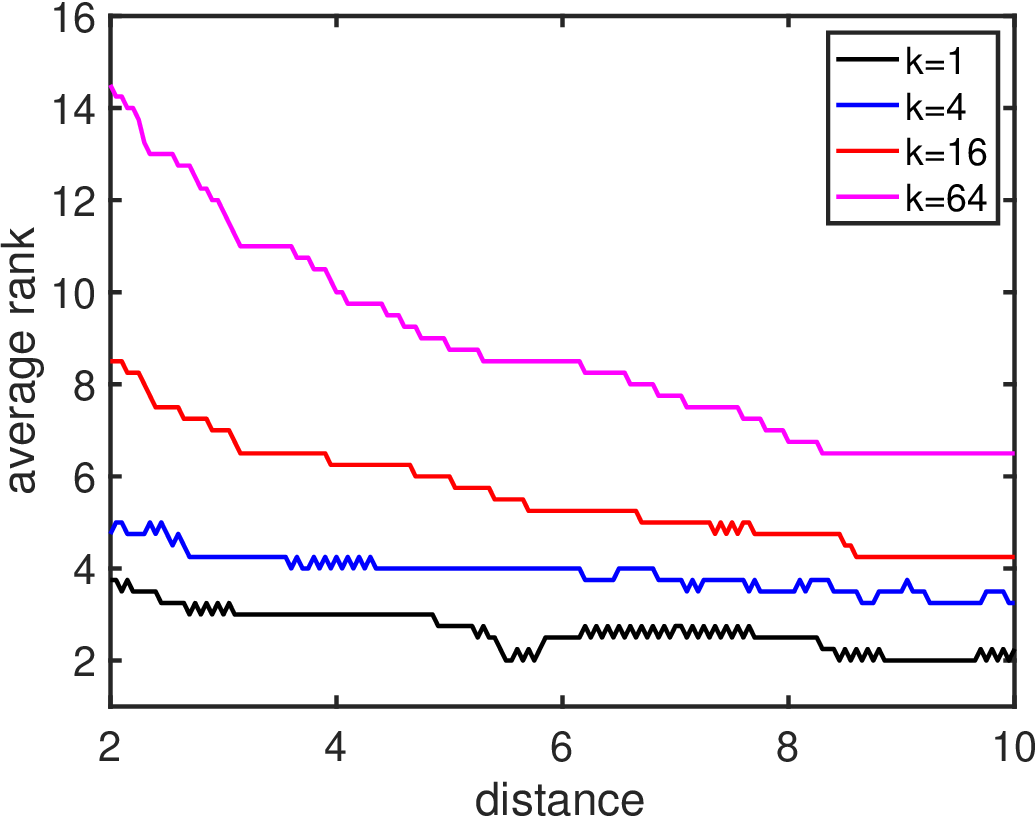}
\includegraphics[width=0.45\textwidth,trim={0.0cm 0.0cm 0.0cm 0.0cm},clip]{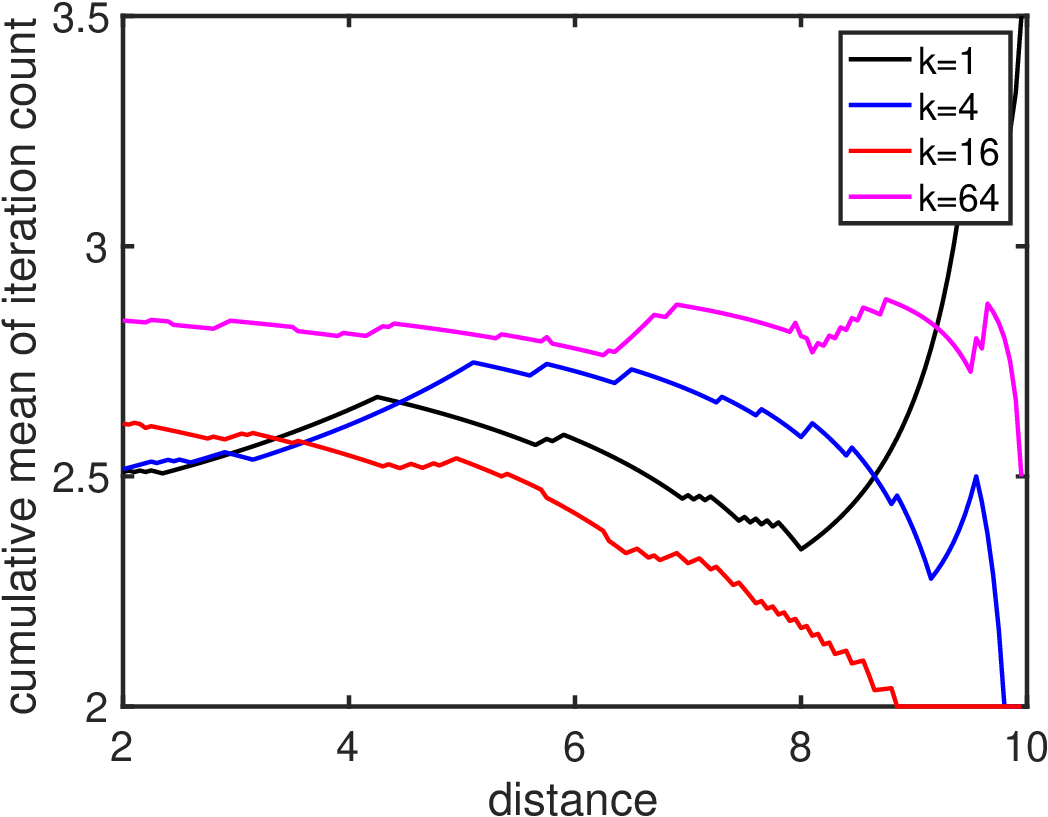}
\includegraphics[width=0.45\textwidth,trim={0.0cm 0.0cm 0.0cm 0.0cm},clip]{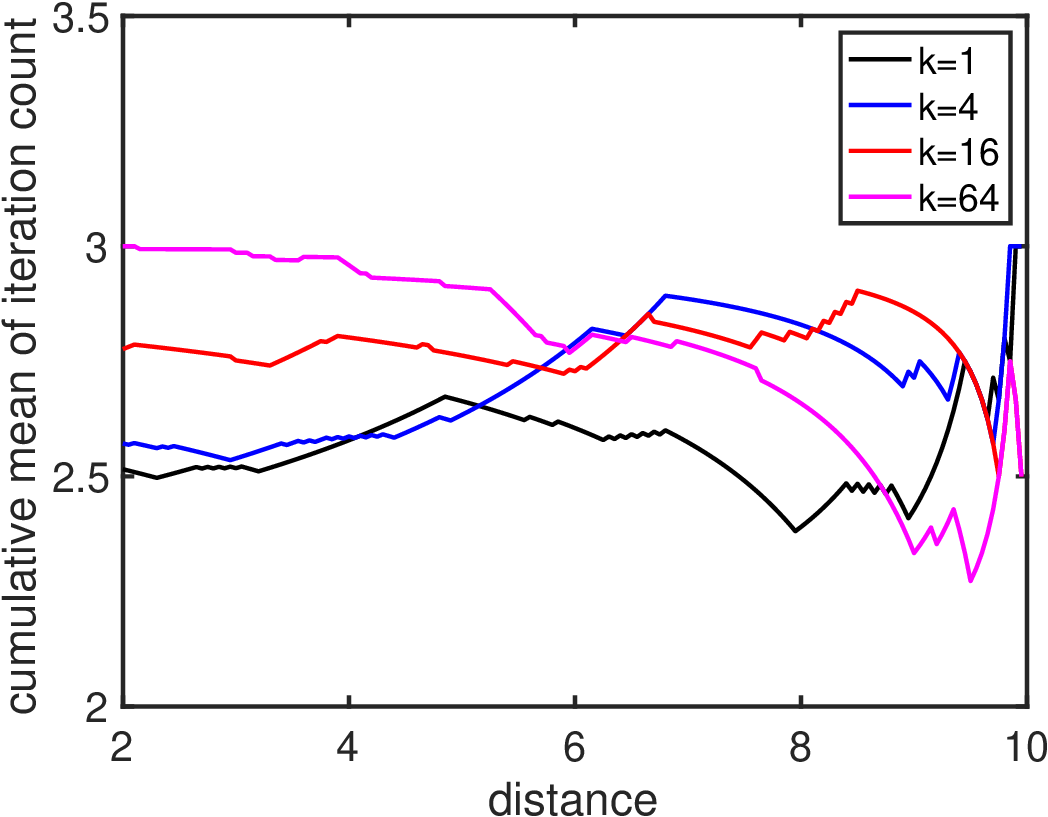}
\caption{Compression of the Hemholtz free space Green's function in three and four dimensions.  Left: three-dimensional results. Right: four-dimensional results.  See the text for details.\label{fig:4D}}
\end{center}
\end{figure}
We expect that C2D is most efficient when it is warm started, i.e. when a good initial guess   is available. This is particularly applicable for parametric tensor approximation.  Instead of testing standard parametric tensors, below we present two interesting applications of C2D for approximating the Green's function solutions to Helmholtz and the solution to Poisson's equation. 

\subsubsection{Compression of the Helmholtz fundamental solution}
In this problem we assume a point source with a time-harmonic forcing is placed in the origin in 3 or 4 dimensions. The free space solution to the Helmholtz equation at wave number $\kappa$ in $d$ dimensions  is then the Green's function 
\begin{equation}
X(x_1,\ldots,x_d) = \frac{i}{4}\left(\frac{\kappa}{2\pi \rho}\right)^{\frac{d}{2}-1} H_{\frac{d}{2}-1}^{(1)}(\kappa \rho). \label{eq:greenHelm}
\end{equation}
 Here $\rho = \sqrt{x_1^2 + x_2^2 + \cdots + x_d^2}$ is the radial distance from the origin and $H_{\frac{d}{2}-1}^{(1)}$ is the Hankel function. It is expected, \cite{parabolic_sep,no_low_rank_Engq}, that the Green's function is easier to compress when the wavenumber is low and when the distance to the source is far away. 

We then try to compress the real part of the Green's function (\ref{eq:greenHelm}) sampled in a unit cube whose corner is placed in the coordinate $(x_1^\ast,0,0)$  in 3D and  $(x_1^\ast,0,0,0)$ in 4D. We sample the function using 100 equidistant points in each dimension. The cube is first placed at $x_1^\ast = 10$ and is then moved towards the origin in steps of 0.05 until its corner closest to the origin is at $x_1^\ast = 2$. For the initial application of C2D we start from a random rank 1 Tucker tensor and for the remaining applications of C2D we use the previous solution at a neighboring $x_1^*$ value as initial guess for warm-start. We repeat the process for three and four dimensions and in both cases we take $\kappa = 1, 4, 16 , 64$, making the function increasingly difficult to compress. In all cases we request a tolerance of $10^{-4}$.

The results of the experiments are displayed in Figure \ref{fig:4D}. To the left we display the results in 3D and to the right we display the results in 4D. The top figures display the average ranks as a function of the distance (in this case $x_1^\ast$) for the four different $\kappa$ (denoted by $k$ in the figure). In the bottom figures we display the cumulative average number of iterations as a function of the distance. 
As expected the Green's function is easier to compress when the wavenumber is low and when the distance to the source is far away. It is encouraging to see that the number of iterations appears to be very manageable, in general between 2 and 3. 

\begin{figure}[h!]
\begin{center}
\includegraphics[width=0.5\textwidth,trim={0.0cm 0.0cm 0.0cm 0.0cm},clip]{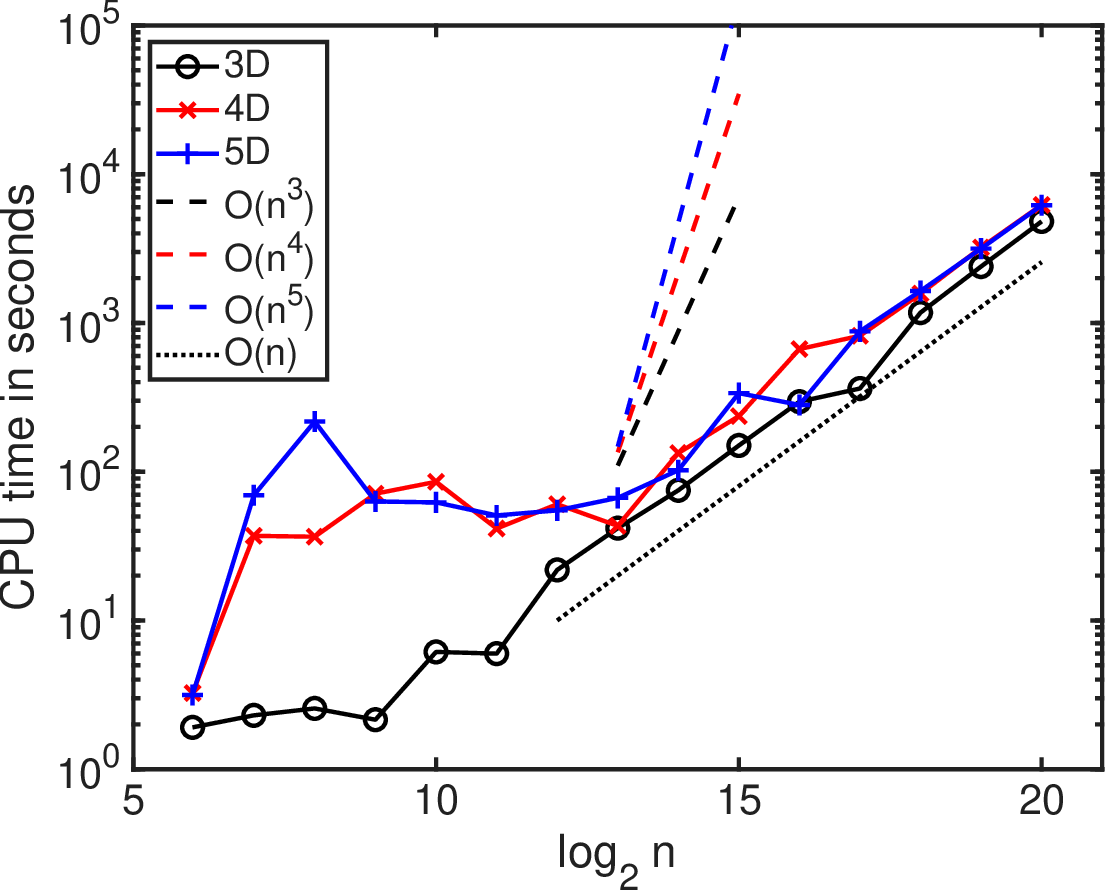}
\caption{Timing results comparing the fast Cross$^2$-DEIM Poisson solver, see the text for details. \label{fig:poisson3D}}
\end{center}
\end{figure}

\subsubsection{A sublinear fast direct Poisson solver using Cross$^2$-DEIM} \label{sec:PoissonDFT}
In this example we show how C2D can be used to construct a very fast direct solver for the Poisson equation in $d$ dimensions 
\begin{equation} \label{eq:poisson}
\Delta v(x_1,\ldots,x_d) = f(x_1,\ldots,x_d), \ \ (x_1,\ldots,x_d) \in [-1,1]^d,
\end{equation}
with homogenous Dirichlet conditions. Here we chose the number of dimensions to be three, four and five and specialize to Dirichlet boundary conditions but note that extension to any dimension and Neumann or periodic boundary conditions is possible.  

For simplicity we assume that each coordinate is discretized using $n+2$ gridpoints on a grid $y_j = -1 + j h, \ \ h = 2/(n+1), j = 1,\ldots,n$. We approximate the derivatives by the standard centered finite difference stencil corresponding to the matrix 

\begin{equation}
\frac{-1}{h^2} T_n =  \frac{-1}{h^2}\begin{bmatrix}
 \phantom{-} 2 & -1 &   \phantom{-} 0 & \cdots &  \phantom{-} 0 & \phantom{-} 0 \\
-1 &  \phantom{-} 2 & -1 &  \phantom{-} 0 & \cdots &  \phantom{-} 0 \\
\vdots & \vdots & \vdots & \vdots & \vdots \\
\phantom{-} 0 &  \phantom{-}  0 & \cdots & -1 &  \phantom{-} 2 & -1 &  \\
\phantom{-} 0 & \phantom{-} & \cdots & 0  & -1 &  \phantom{-} 2
\end{bmatrix}.
\end{equation}
It is well know, \cite{Demmel:1997pt}, that the eigenvalues, $\lambda_j$ and the eigenvectors $z_j(k)$ to $T_n$ are 
\[
\lambda_j = 2(1-\cos(\frac{\pi j}{n+1}), \ \ z_j(k) = \sqrt{\frac{2}{n+1}} \sin \left(\frac{j k \pi}{n+1} \right),
\]
so that $T_n$ has the eigen-decomposition $T_n = Z \Lambda Z^T$. As described e.g. in the textbook by Demmel, \cite{Demmel:1997pt}, the application of the eigenvector matrix $Z$ can be done by a fast sine transform. Then, for a vector $w\in \mathbb{R}^n$ we have 
\[
Z w = \sqrt{\frac{2}{n+1}} \verb+dst+(w),
\]   
where $\verb+dst+(\cdot)$ is the application of the fast sine transform. 

We explain our approach in three dimensions, while noting the extensions to higher dimensions are trivial.  Suppose that $F = G_{F} \times_1 U_1 \times_2 U_2  \times_3 U_3$ is a Tucker approximation   to the right hand side of (\ref{eq:poisson}) then, following the full rank ``FFT solver'' in \cite{Demmel:1997pt},  we can define 
\[
\hat{F} = G_{F} \times_1 \sqrt{\frac{2}{n+1}} \verb+dst+(U_1) \times_2 \sqrt{\frac{2}{n+1}} \verb+dst+(U_2)  \times_3 \sqrt{\frac{2}{n+1}} \verb+dst+(U_3).
\]
Let $\hat{V} \approx G_{V} \times_1 \hat{U}_1 \times_2 \hat{U}_2  \times_3 \hat{U}_3$ be the Cross$^2$-DEIM approximation of the tensor $\hat{V}$ whose elements are given by
\begin{equation}
\label{c2dwv}
\hat{V}(i_1,i_2,i_3) = \frac{-h^2\hat{F}(i_1,i_2,i_3)}{\lambda_{i_1}+ \lambda_{i_2} + \lambda_{i_3}}.
\end{equation}
Then an approximate solution to (\ref{eq:poisson}) is given by 
\[
V = G_{V} \times_1 \sqrt{\frac{2}{n+1}} \verb+dst+(\hat{U}_1) \times_2 \sqrt{\frac{2}{n+1}} \verb+dst+(\hat{U}_2)  \times_3 \sqrt{\frac{2}{n+1}} \verb+dst+(\hat{U}_3).
\]
To compute the numerical solution $V$ we thus need a few fast sine transforms and a single C2D approximation of \eqref{c2dwv}. 

To demonstrate the numerical performance of this solver, we compute approximate solutions using $n = 2^p-1$, $p = 6,7,\ldots,20$, grid points in each direction. The right hand side is chosen to be $f = e^{-36\rho^2}$, $\rho^2 = \sum_{i = 1}^d (x_i-i/100)^2$. For dimensions three to five, we solve, in sequence, for each $n$. For the coarsest mesh we use the factor matrices $\verb+dst+(\hat{U}_i), i=1, \ldots, d$ to warm-start the method, and for the remaining values of $n$, we prolongate the previous solution to the finer grid using linear interpolation, and use the prolongated solution for the warm start. We chose the tolerance so that the rank is approximately 5 in each mode for the three and four dimensional cases and so that it is three for the five dimensional case. When using the prolongated solution for the warm start C2D requires 3-4 iterations to converge. 

We record the CPU times (in seconds) for each $n$ and $d$ on a M2 MacBook Pro with 16GB RAM. We note that the implementation is in Matlab using the package Tensorlab and that our naive implementation evaluates each element needed in C2D without any reuse and that Tensorlab / Matlab does not exploit the multiple cores on the M2 chip for the element-wise evaluation. If the rank is $r$ the evaluation of each element costs $dr^d$ leading to a cost of $\mathcal{O}(d^2nr^{d+1} + dr^{2d})$ for each iteration in C2D. This means that for this example, sampling cost   will dominate the computational cost in C2D.       
The results are displayed in Figure \ref{fig:poisson3D} where we see that the linear scaling in $n$ becomes dominant for $n$ around $2^{12}=4096$, while for coarser mesh size, the CPU time has very weak dependence on $n.$ In the figure the dashed lines indicate optimal complexity $n^d$ for a full rank solution method (e.g. using multigrid). Considering the case $n=2^{12}$ (setting aside that the full rank solution itself could not be stored on a laptop) and assuming an optimistic time to solution per degree of freedom of 1 nano second it would take around a minute, 78 hours, and 36 years to find the solution for the three, four and five dimensional problems, respectively.

\subsection{Numerical experiments with Tucker-AA}
In this section we present examples with the Cross$^2$-DEIM enabled Tucker-AA.
\subsubsection{Bratu problem}
In this example we use Tucker-AA to solve the non-linear Bratu problem in three dimensions   
\[
\Delta v(x_1,x_2,x_3) + \lambda e^{v} = 0, \ \ (x_1,x_2,x_3) \in [0,1]^3,
\]
with $\lambda=1$ and homogeneous Dirichlet boundary conditions. To find an approximate solution we discretize this equation using standard second order finite difference approximations for the spatial derivatives. We use a uniform grid with grid size $h = 1/(n+1)$ in all dimensions. We have taken $n=256$ in this example. Given an approximation $X(i_1,i_2,i_3) \approx v(i_1h, i_2h ,i_3h), \, i_1,i_2,i_3 = 0,\ldots,n+1,$ this results in a function $H_{\rm B}(i_1,i_2,i_3;X)$ describing the finite difference approximation  
\begin{gather}
\begin{split}
H_{\rm B}(i_1,i_2,i_3;X) &= \frac{1}{h^2}\left(X(i_1+1,i_2,i_3) - 2X(i_1,i_2,i_3) + X(i_1-1,i_2,i_3)  \right) \\
&+ \frac{1}{h^2} \left(X(i_1,i_2+1,i_3) - 2X(i_1,i_2,i_3) + X(i_1,i_2-1,i_3)  \right) \\
&+ \frac{1}{h^2} \left(X(i_1,i_2,i_3+1) - 2X(i_1,i_2,i_3) + X(i_1,i_2,i_3-1)  \right) + \lambda e^{X(i_1,i_2,i_3)}.
\end{split}
\end{gather}
Near the boundaries some of the terms in this expression will be set to zero to account for the homogeneous Dirichlet boundary conditions. 

The fixed point function in Tucker-AA, $H(i_1,i_2,i_3)$,  is obtained by applying the preconditioned Richardson iteration. We have 
\[
X_{k+1}(i_1,i_2,i_3) = H(i_1,i_2,i_3;X_{k},\alpha) \equiv  X_k(i_1,i_2,i_3) + \alpha M(H_{\rm B}(i_1,i_2,i_3;X_{k})).
\] 
Here $M$ represents the preconditioner, which, when used, we take to be the fast Poisson solver described in Section \ref{sec:PoissonDFT} and use the current solution $X_k$ for the warm start. For the preconditioner we use a C2D tolerance that is ten times smaller than the current Tucker-AA rounding tolerance, i.e. ${\rm TOL} = 0.1 \epsilon_H$.

\begin{figure}[]
\begin{center}
\includegraphics[width=0.60\textwidth,trim={0.0cm 0.0cm 0.0cm 0.0cm},clip]{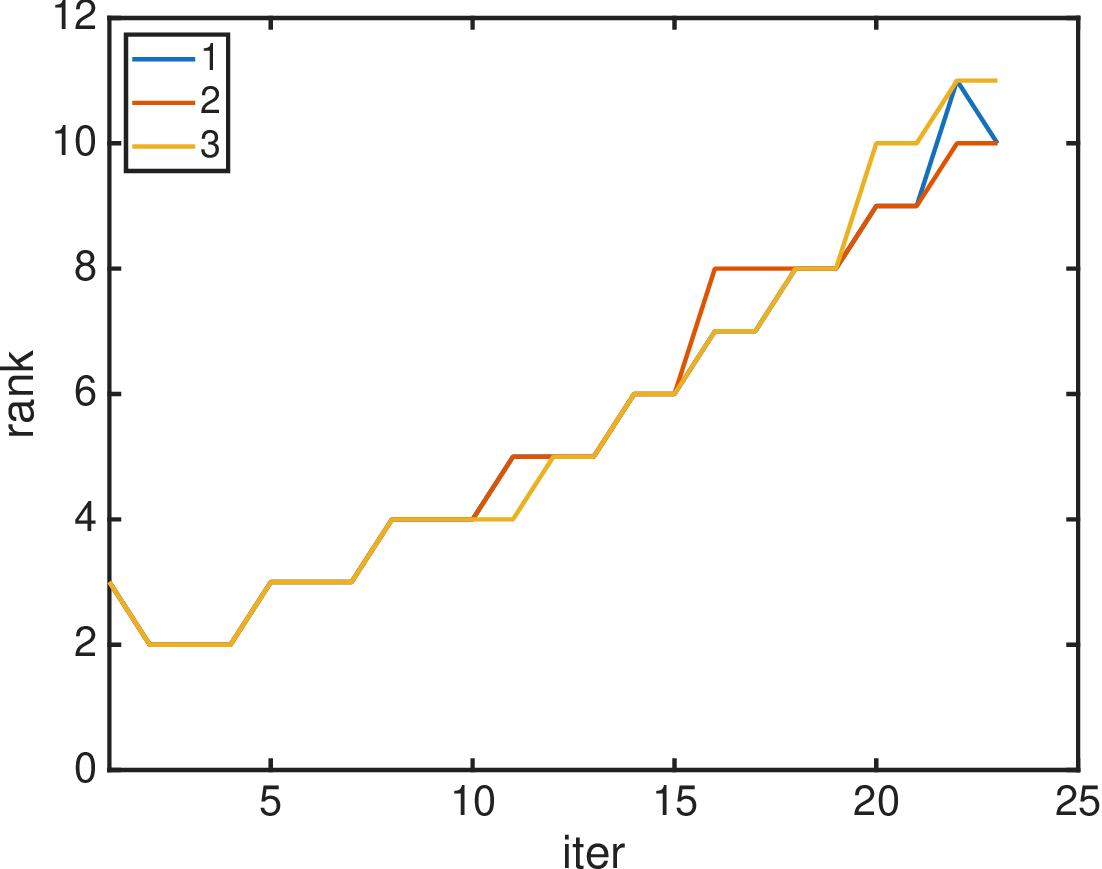} \ \ \ \ 
\includegraphics[width=0.335\textwidth,trim={0.0cm 0.0cm 0.0cm 0.0cm},clip]{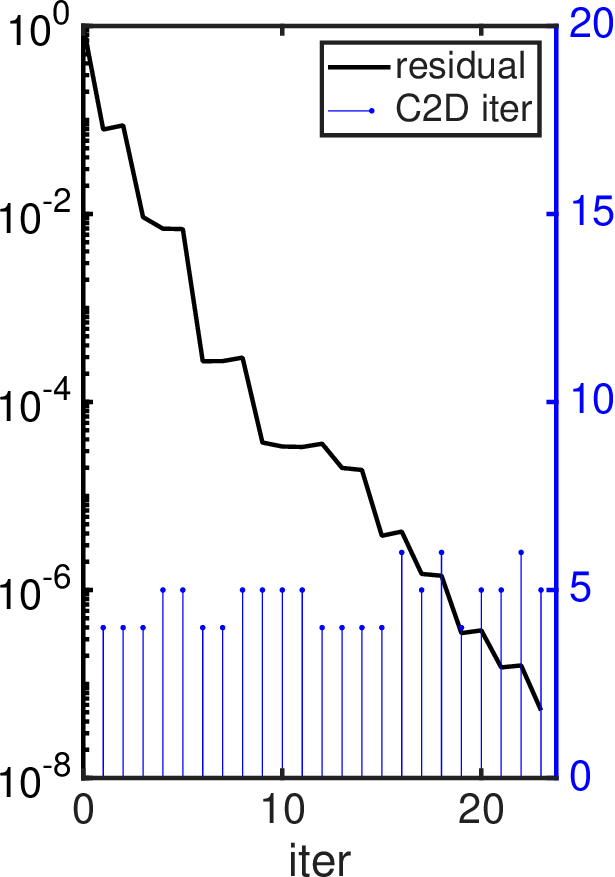} 
\caption{The figure displays results from the preconditioned Bratu problem. To the left is the ranks of the intermediate Tucker-AA iterates as a function of the iteration count in modes 1, 2, 3. To the right is the residual and the number of C2D iterations in the approximation of $H(i_1,i_2,i_3;X_{k},\alpha)$ as functions of the iteration count. \label{fig:bratu_pre}}
\end{center}
\end{figure}

We first test the Tucker-AA with the fast Poisson preconditioner and with no scheduling. We set ${\rm TOL}=10^{-7} \rho_0, \hat{m}=4, $ and   $\alpha = 0.1$. As initial data we use the one dimensional solution $2 \log \left( \frac{\cosh (\theta) }{\cosh (\theta (1-2x))} \right)$, where $\theta$ solves $\cos (\theta) = 4 \frac{\theta }{\sqrt{\lambda}}$.

To the left in Figure \ref{fig:bratu_pre} we display the ranks of the three unfoldings of the core tensor as a function of the number of iterations. As can be seen the ranks increase monotonically and the solution converges to desired tolerance in a small number of iterations.  To the right in Figure \ref{fig:bratu_pre} we display the normalized residual $\rho_k/\rho_0$ as a function of the iteration number along with the number of C2D iterations (note that the iteration count is displayed on the right axis) needed to approximate $H(i_1,i_2,i_3;X_{k},\alpha)$. We notice rapid and monotone decay of the residual, and the C2D iteration numbers are around 5.

\begin{figure}[]
\begin{center}
\includegraphics[width=0.32\textwidth,trim={0.0cm 0.0cm 0.0cm 0.0cm},clip]{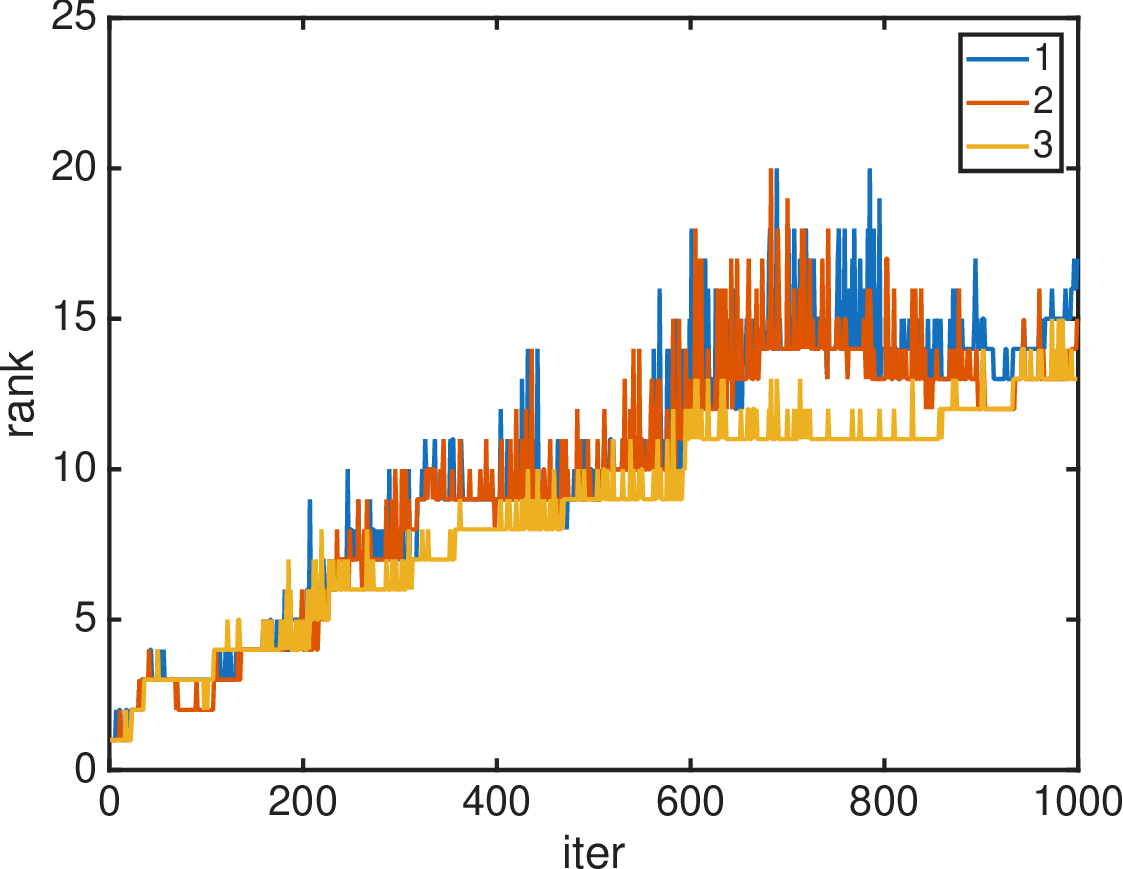} 
\includegraphics[width=0.32\textwidth,trim={0.0cm 0.0cm 0.0cm 0.0cm},clip]{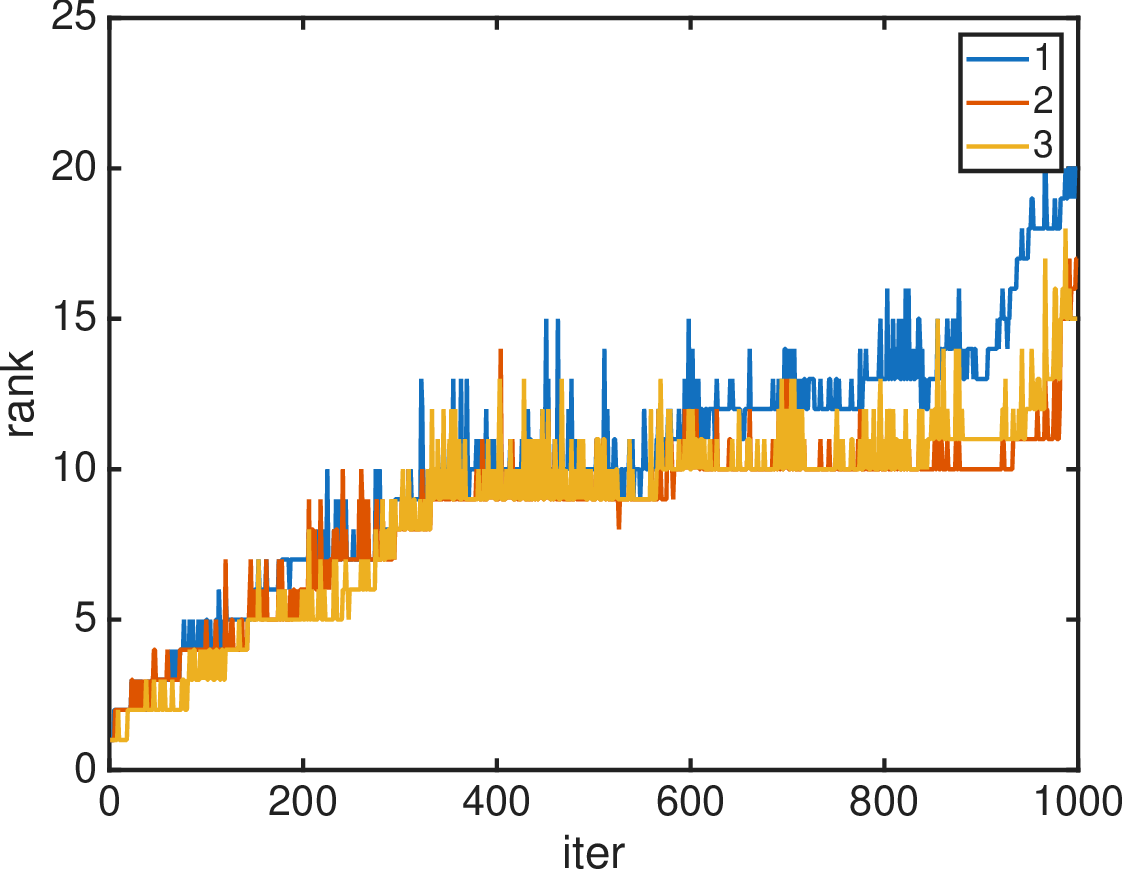} 
\includegraphics[width=0.32\textwidth,trim={0.0cm 0.0cm 0.0cm 0.0cm},clip]{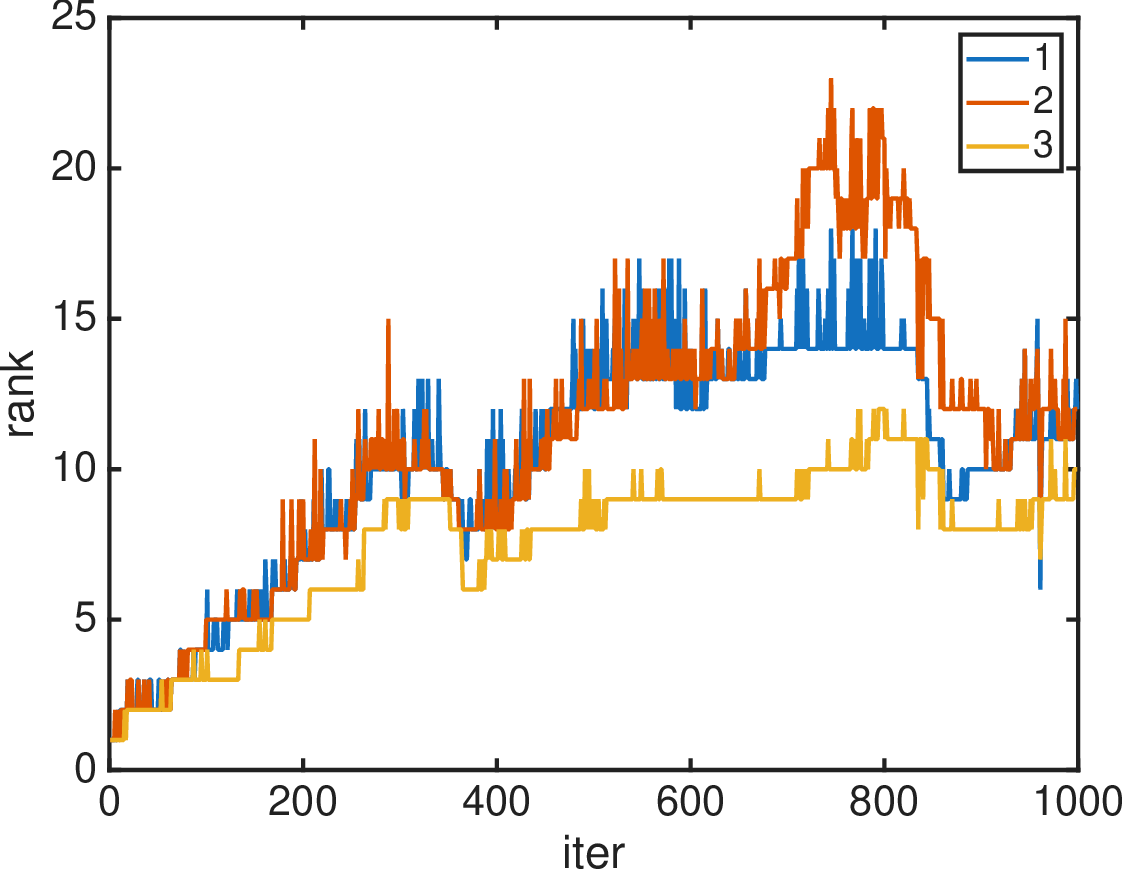} 
\includegraphics[width=0.32\textwidth,trim={0.0cm 0.0cm 0.0cm 0.0cm},clip]{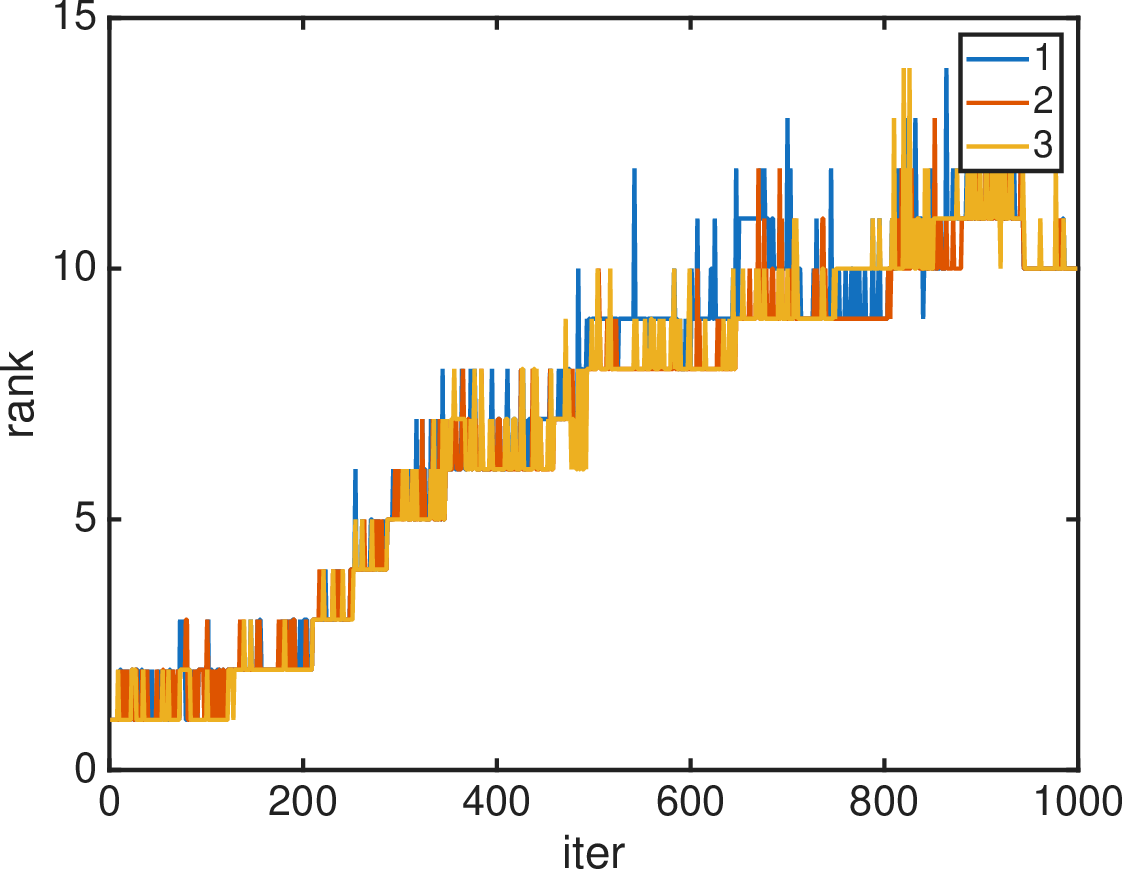} 
\includegraphics[width=0.32\textwidth,trim={0.0cm 0.0cm 0.0cm 0.0cm},clip]{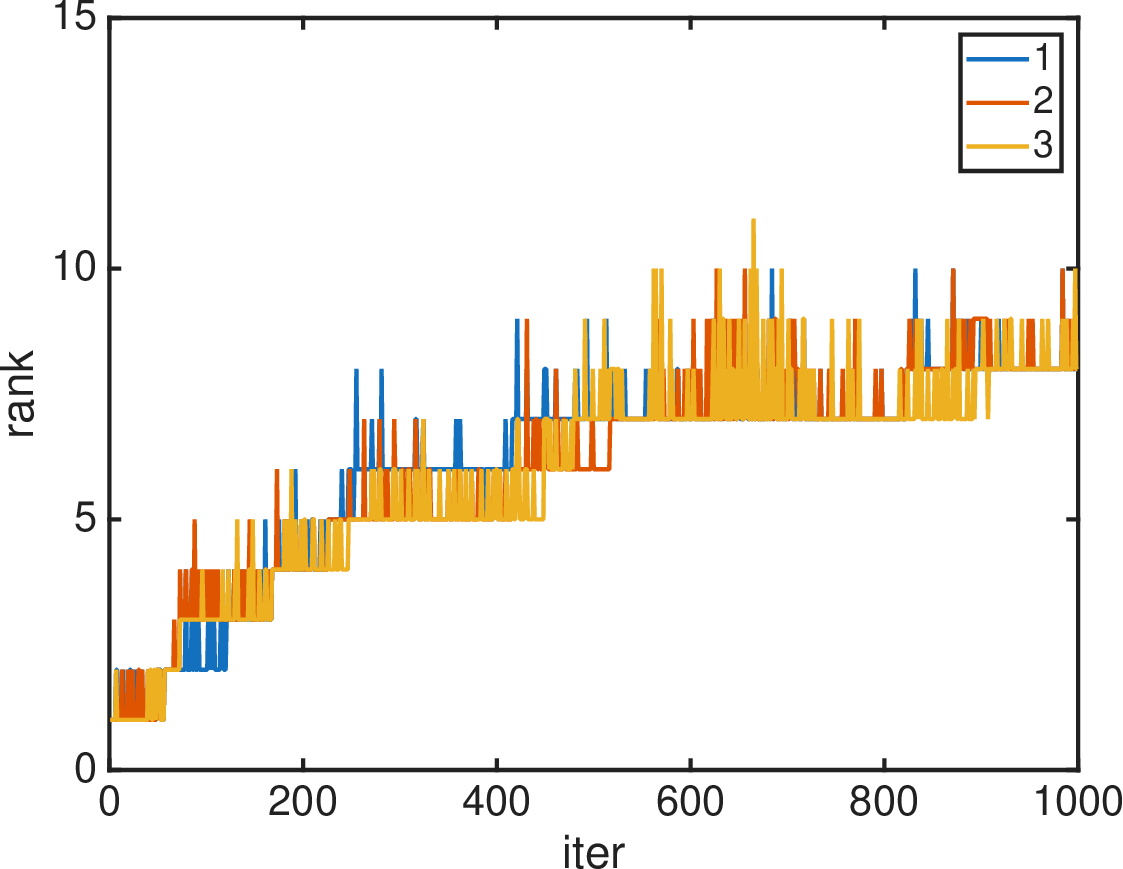} 
\includegraphics[width=0.32\textwidth,trim={0.0cm 0.0cm 0.0cm 0.0cm},clip]{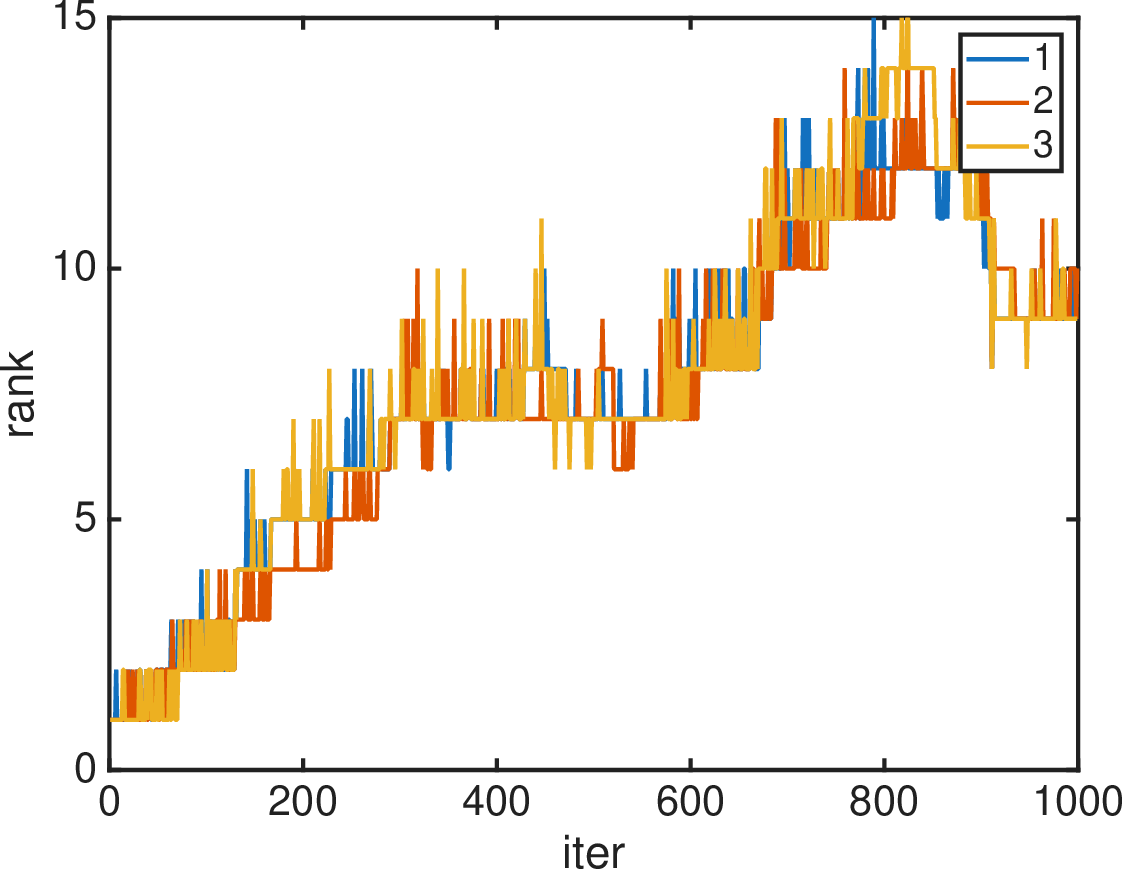} 
\caption{The figure displays the ranks of the iterates $X_{k}$ for the Bratu problem with no preconditioning for modes 1, 2, 3. The top row is for $\theta =0.5$ and the bottom row is for $\theta = 0.9$. From left to right we cap the number of C2D iterations to 2, 4, and 10.   Note that the scale of the vertical axis is different for the two rows. \label{fig:bratu_rank}}
\end{center}
\end{figure}

For the same problem (using the same initial guess) we use Tucker-AA with no preconditioning and vary the scheduling parameter $\theta$, the window size $\hat{m}$ and the maximum number of C2D iterations. In all computations we monitor the ranks throughout the iteration and at the final iteration, the residual and the number of C2D iterations. We stop the Tucker-AA iteration after 1000 iterations. 

\begin{figure}[]
\begin{center}
\includegraphics[width=0.16\textwidth,trim={0.0cm 0.0cm 0.0cm 0.0cm},clip]{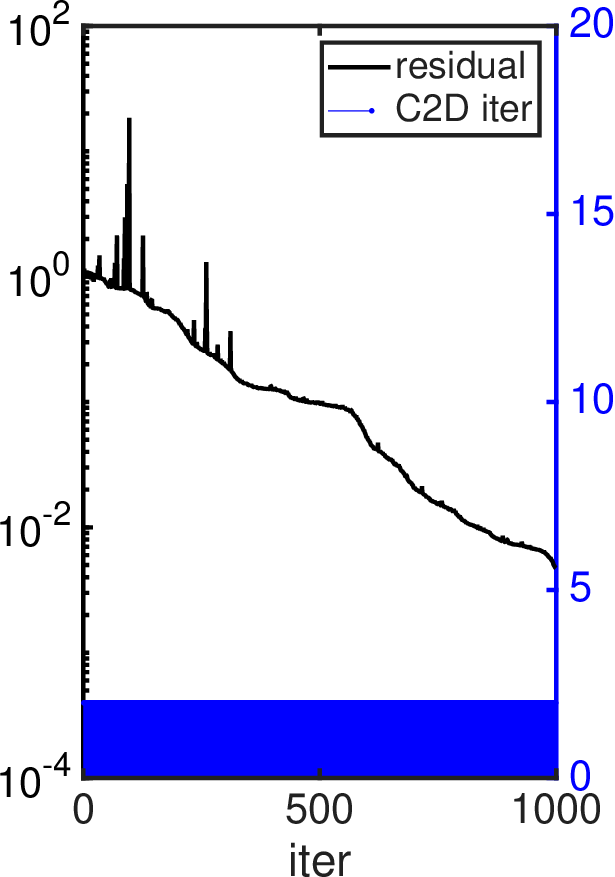} 
\includegraphics[width=0.16\textwidth,trim={0.0cm 0.0cm 0.0cm 0.0cm},clip]{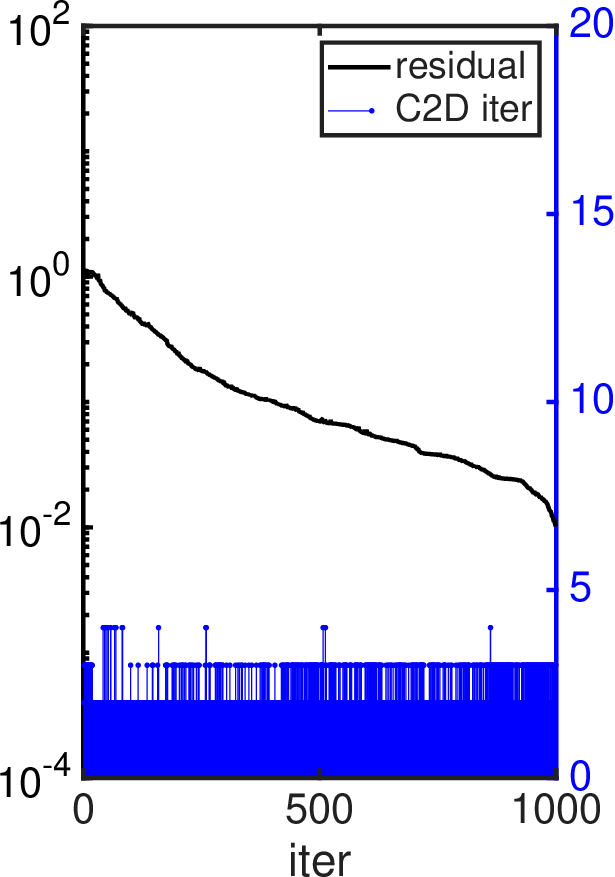} 
\includegraphics[width=0.16\textwidth,trim={0.0cm 0.0cm 0.0cm 0.0cm},clip]{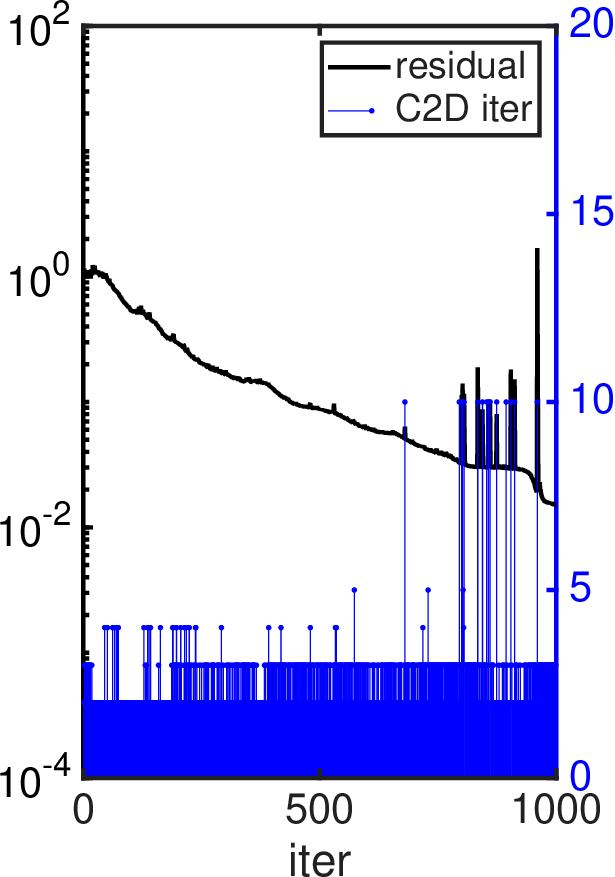} 
\includegraphics[width=0.16\textwidth,trim={0.0cm 0.0cm 0.0cm 0.0cm},clip]{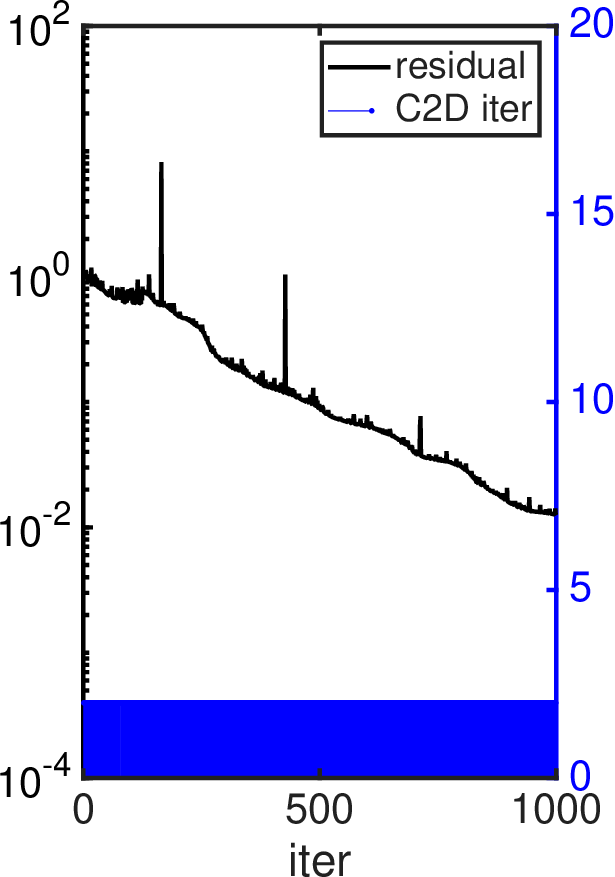} 
\includegraphics[width=0.16\textwidth,trim={0.0cm 0.0cm 0.0cm 0.0cm},clip]{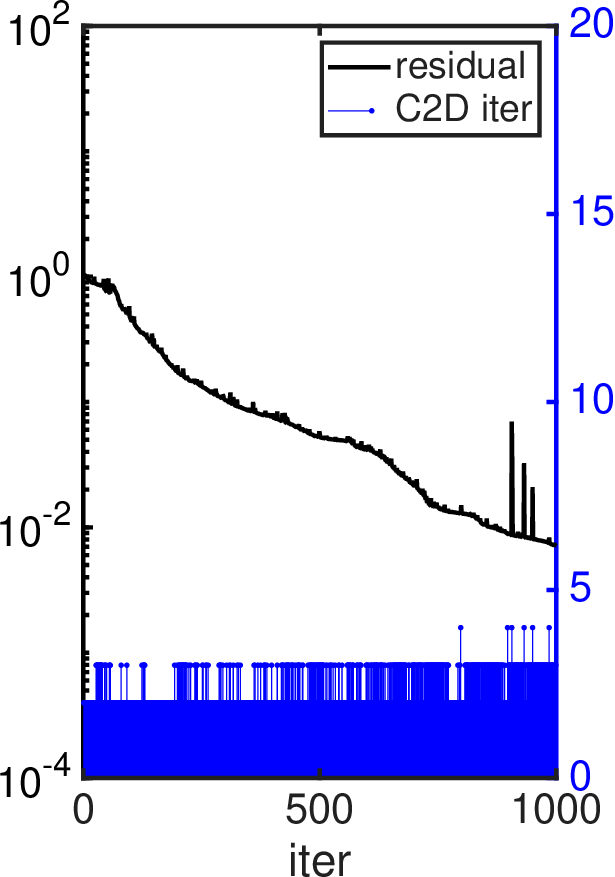} 
\includegraphics[width=0.16\textwidth,trim={0.0cm 0.0cm 0.0cm 0.0cm},clip]{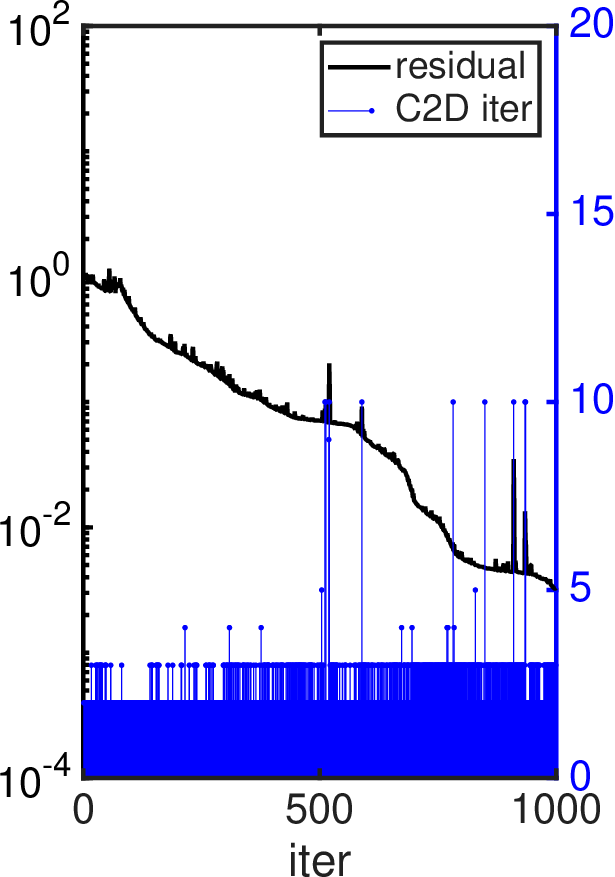} 
\caption{The figure displays the residuals (left vertical axis) and the number of C2D iterations (right axis) of the iterates $X_{k}$ for the Bratu problem with no preconditioning. The first three figures are for $\theta =0.5$ and three next are for $\theta = 0.9$. In the sets of three (from left to right) we cap the number of C2D iterations to 2, 4, and 10. \label{fig:bratu_resid}}
\end{center}
\end{figure}

We first fix $\hat{m} = 4$ and consider $\theta=0.5$  and $\theta=0.9$. For all cases we set $\alpha = 0.1 h^2$. For each $\theta$ we also consider C2D with a maximum of 2, 4 and 10 iterations in an effort to characterize the how this parameter affects the convergence of Tucker-AA and the intermediate rank growth of $X_{k}$. Figure \ref{fig:bratu_rank} displays the ranks as a function of the iteration and Figure \ref{fig:bratu_resid} displays the residual and number of C2D iterations for $H(i_1,i_2,i_3;X_{k},\alpha)$. As can be seen taking $\theta$ closer to one helps with controlling the rank as does restricting the number of C2D iterations to 2 or 4. It can also be seen that for the iterates where C2D requires a large number of iterations the residual often has a sudden large increase. However, it does not appear as if these large residuals affect the overall convergence dramatically. Since the cost of an iteration scales with the rank the overall best parameter choice (for this example) is $\theta =0.9$ in combinations with a cap of 4 on the number of C2D iterations. 

\begin{figure}[]
\begin{center}
\includegraphics[width=0.32\textwidth,trim={0.0cm 0.0cm 0.0cm 0.0cm},clip]{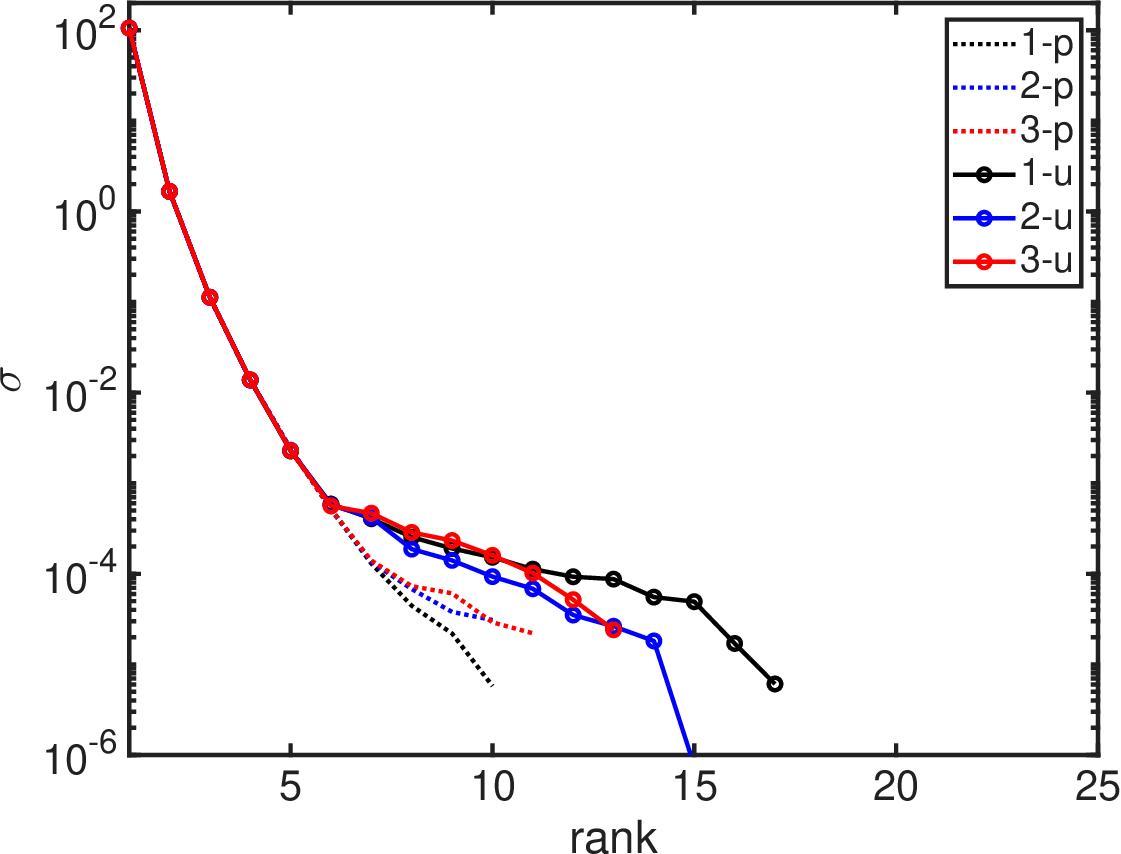} 
\includegraphics[width=0.32\textwidth,trim={0.0cm 0.0cm 0.0cm 0.0cm},clip]{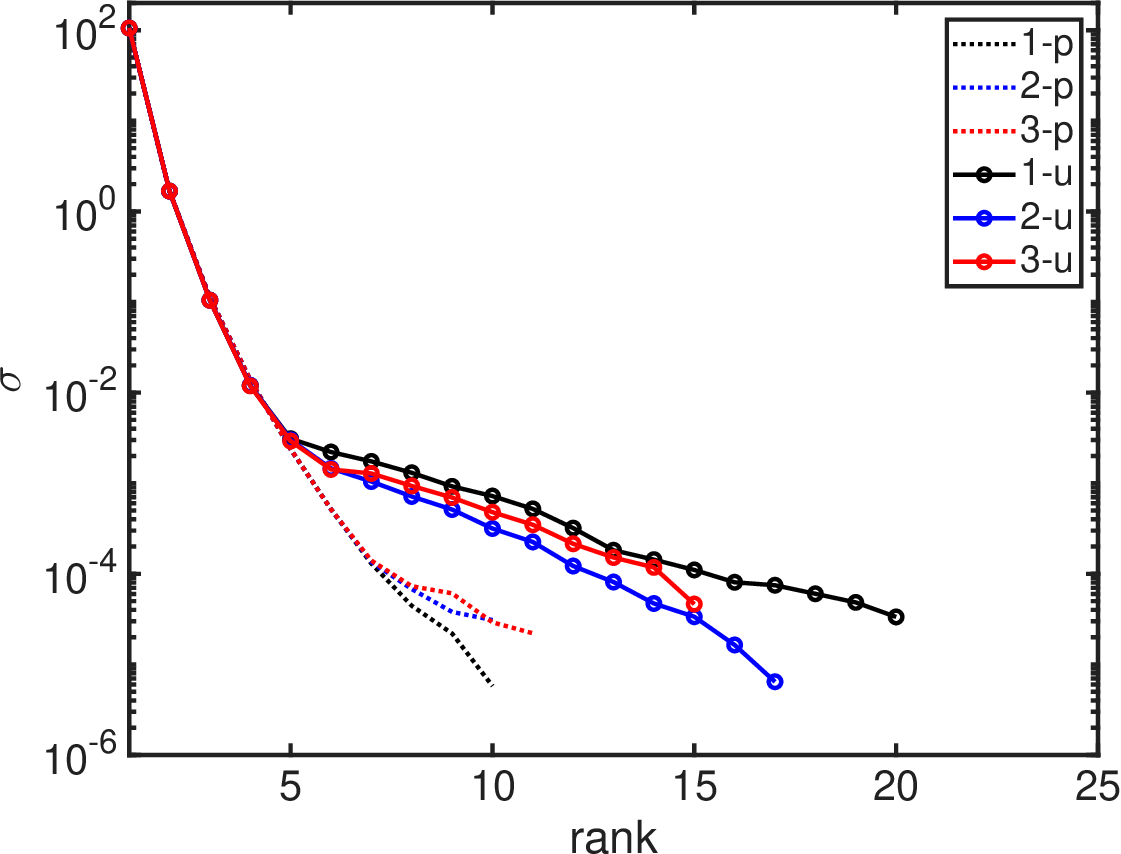} 
\includegraphics[width=0.32\textwidth,trim={0.0cm 0.0cm 0.0cm 0.0cm},clip]{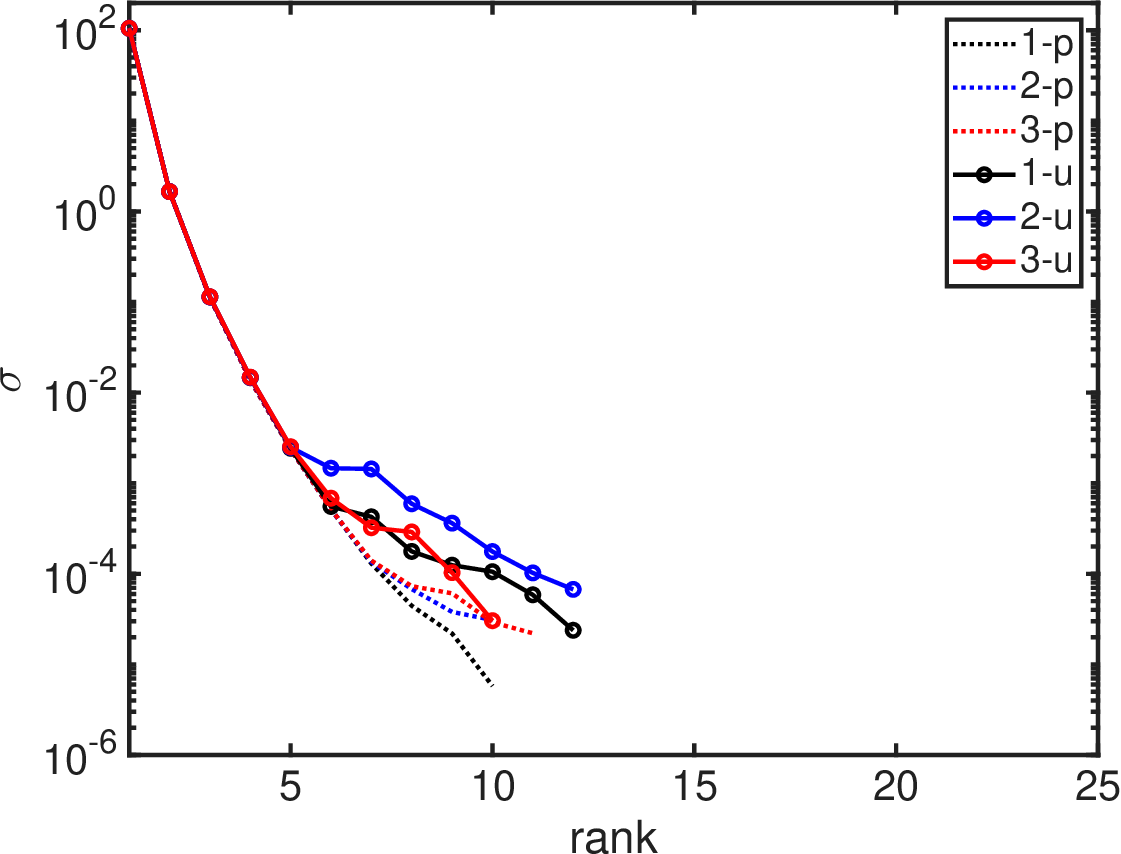} 
\includegraphics[width=0.32\textwidth,trim={0.0cm 0.0cm 0.0cm 0.0cm},clip]{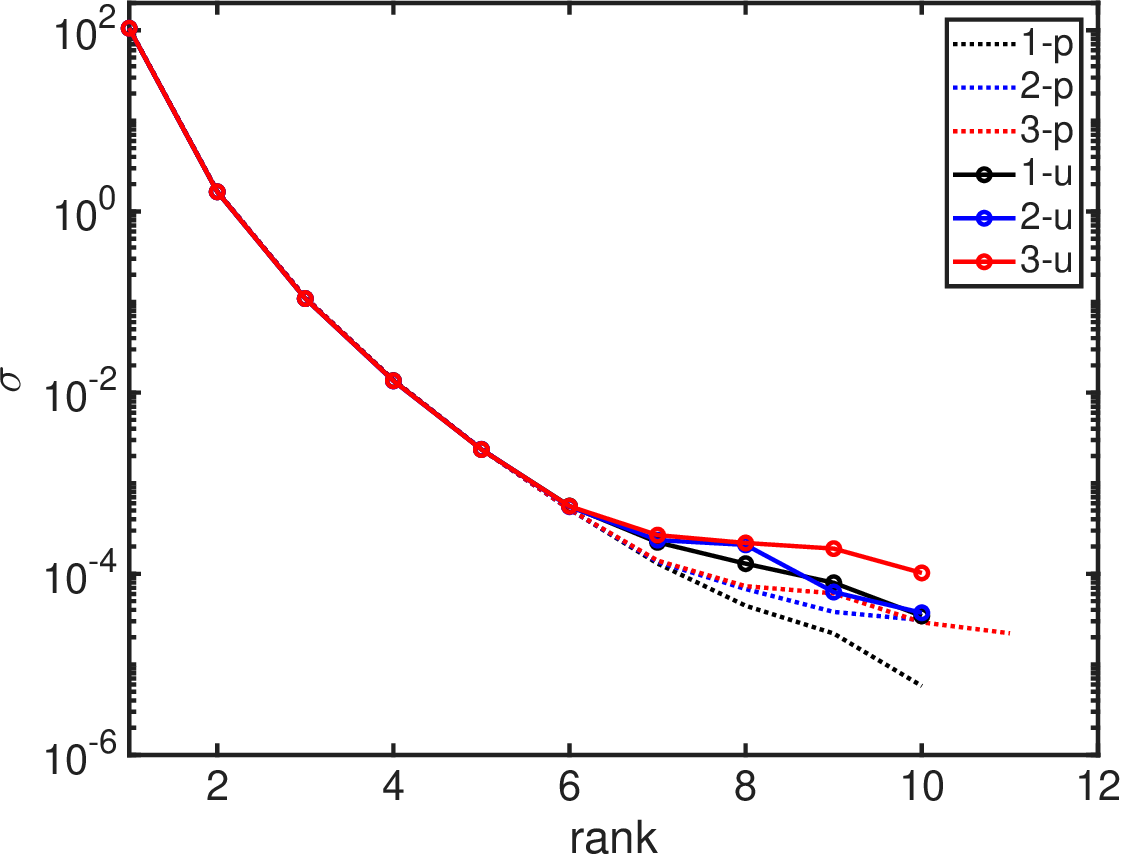} 
\includegraphics[width=0.32\textwidth,trim={0.0cm 0.0cm 0.0cm 0.0cm},clip]{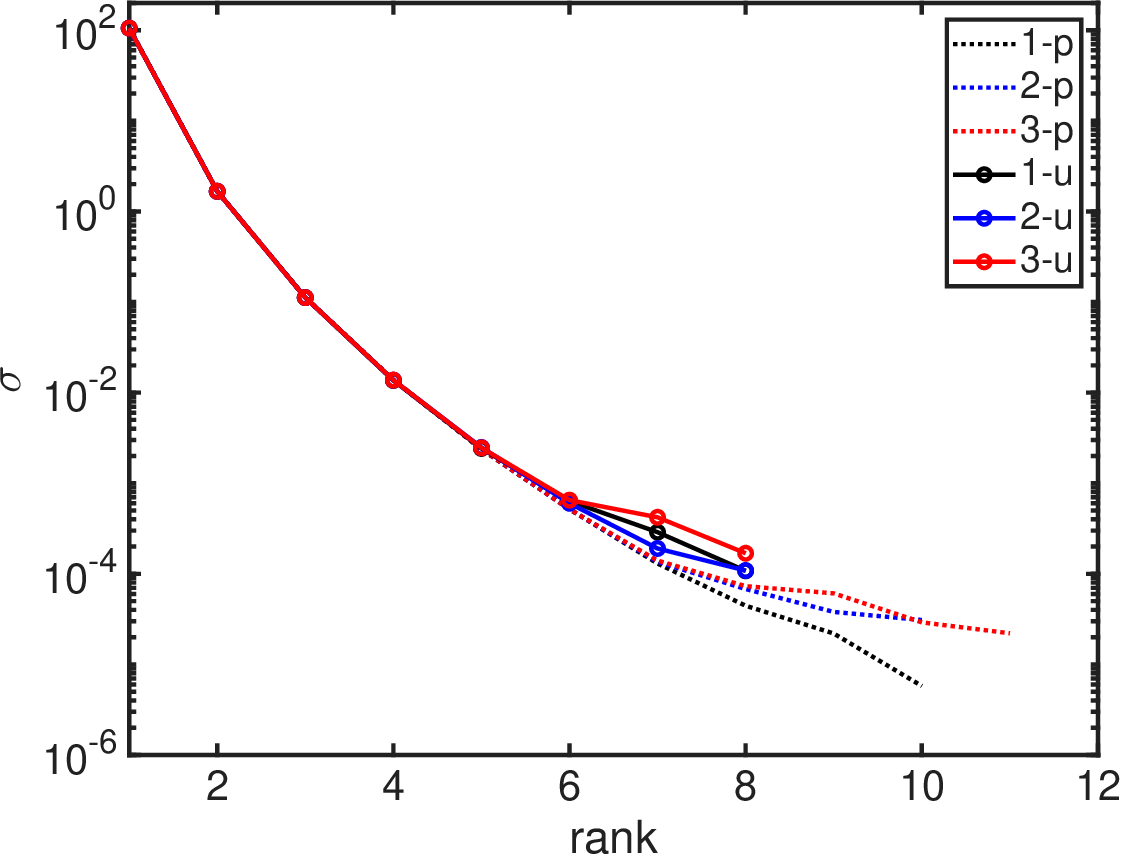} 
\includegraphics[width=0.32\textwidth,trim={0.0cm 0.0cm 0.0cm 0.0cm},clip]{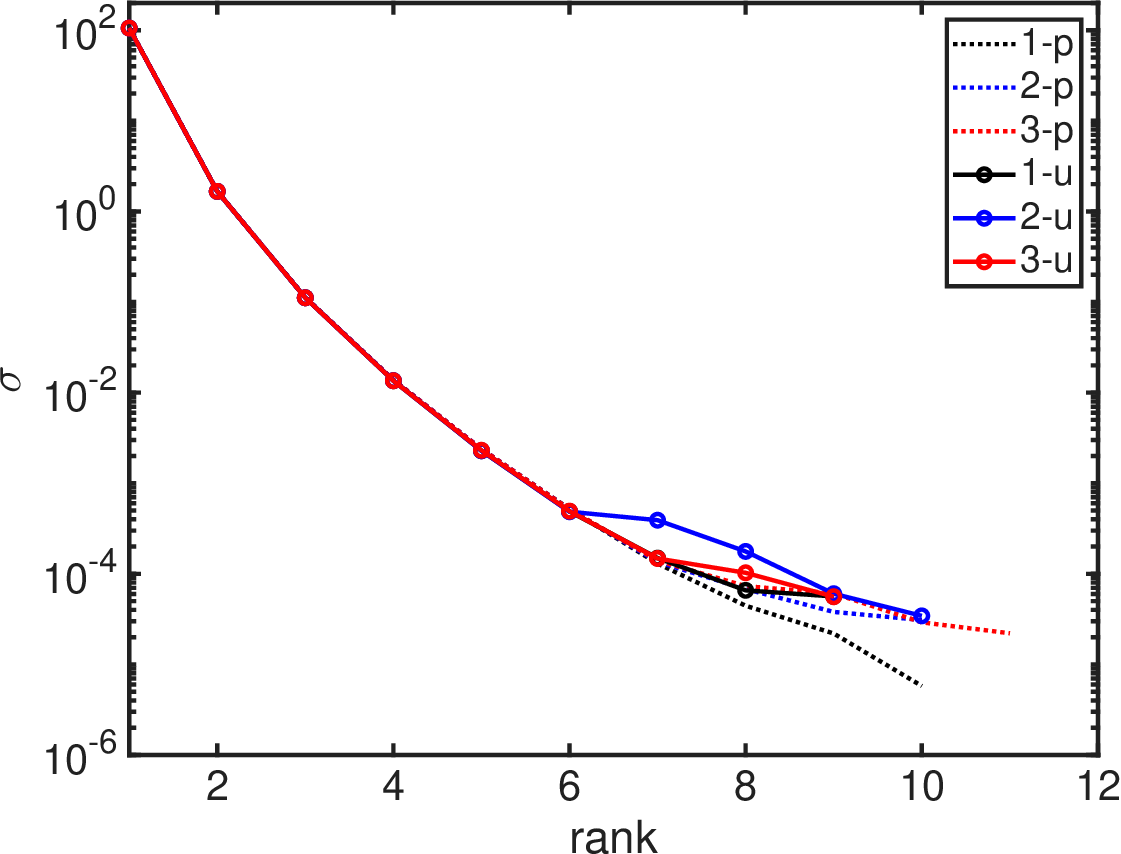} 
\caption{Displayed are the singular value distribution for the unfolding of the core in all directions and for the final iterate together with the singular values for the final   solution obtained with the un-conditioned Tucker-AA. (The dots are obtained with the preconditioned Tucker-AA as reference solution.) The ordering of the graphs is the same as in Figure \ref{fig:bratu_rank}. \label{fig:bratu_rank_final}}
\end{center}
\end{figure}

\begin{figure}[]
\begin{center}
\includegraphics[width=0.32\textwidth,trim={0.0cm 0.0cm 0.0cm 0.0cm},clip]{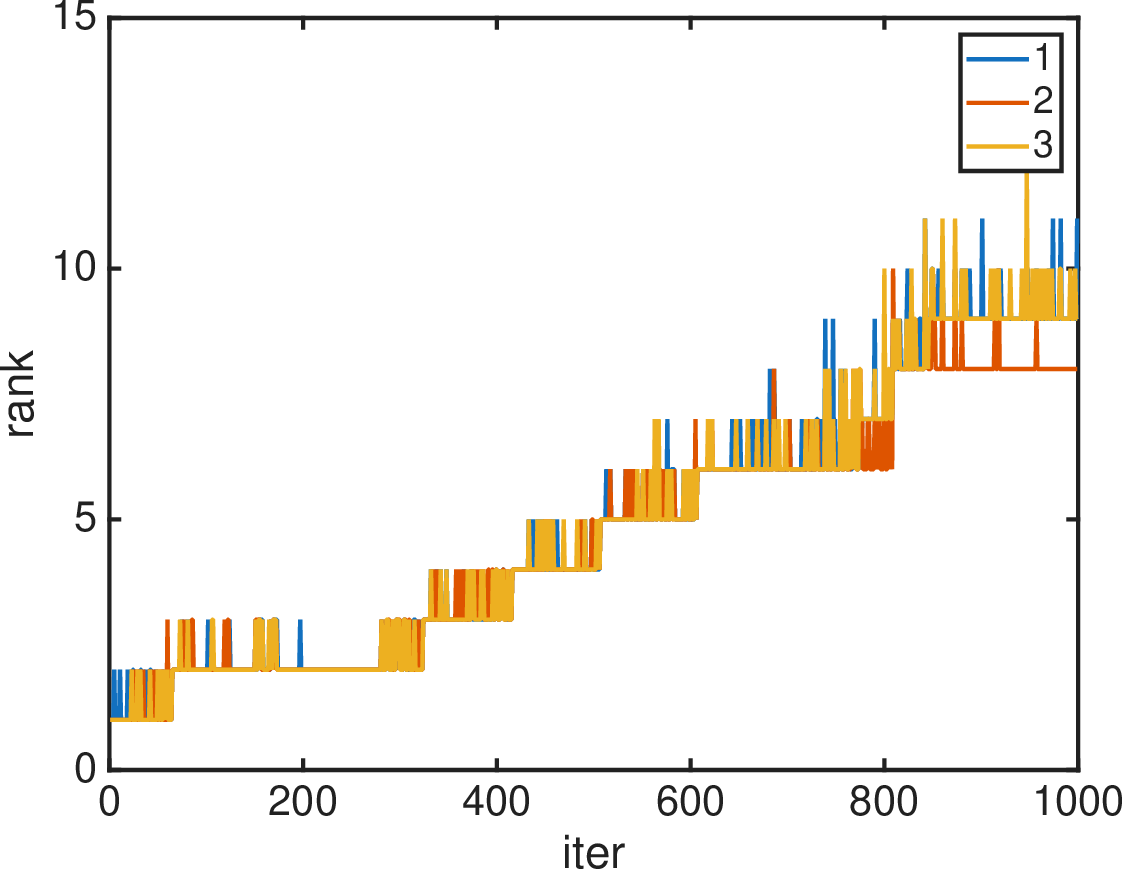} 
\includegraphics[width=0.32\textwidth,trim={0.0cm 0.0cm 0.0cm 0.0cm},clip]{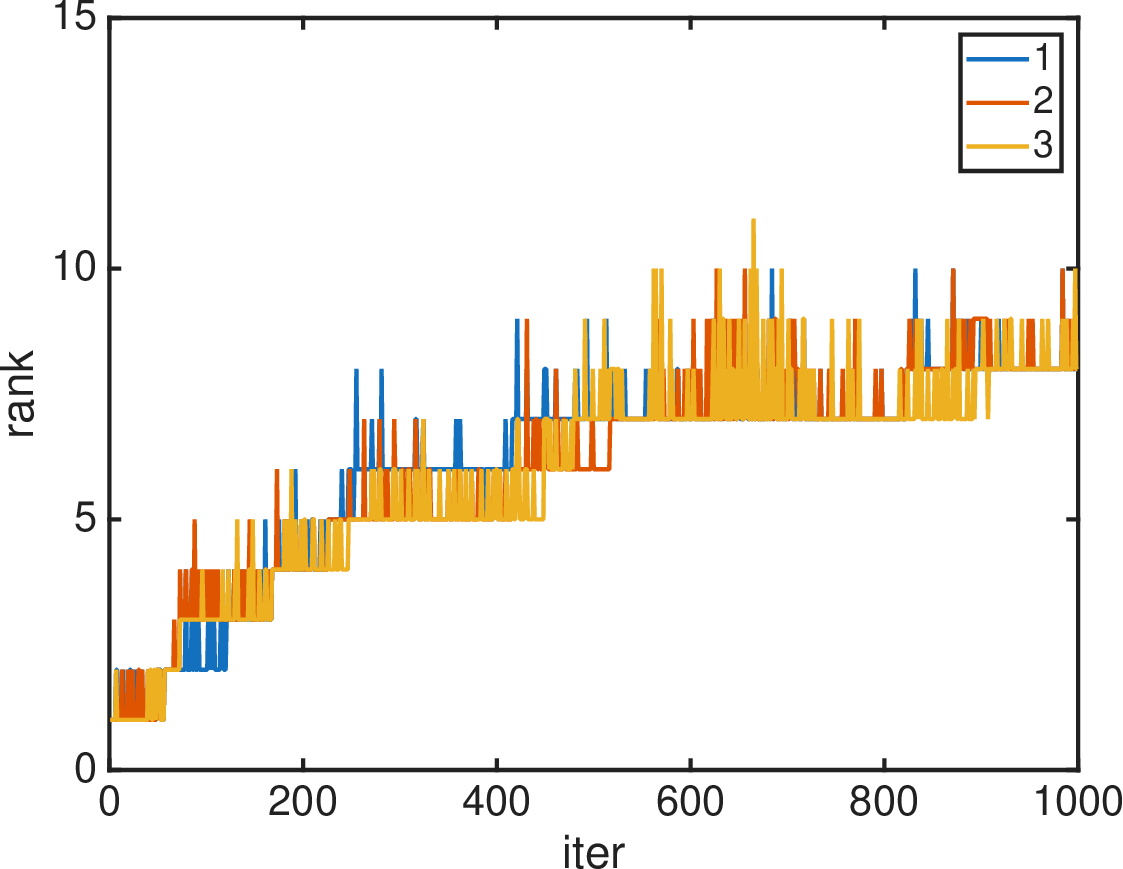} 
\includegraphics[width=0.32\textwidth,trim={0.0cm 0.0cm 0.0cm 0.0cm},clip]{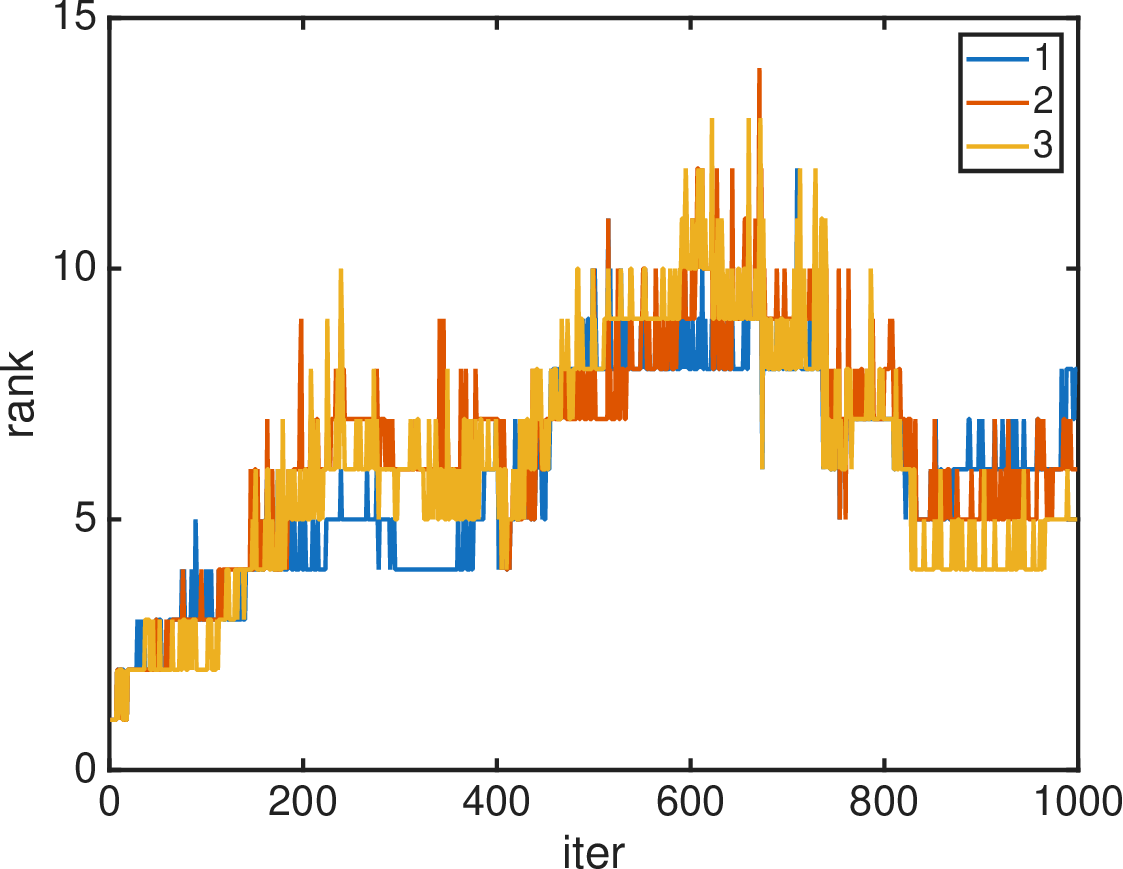} 
\caption{The figure displays the ranks of the iterates $X_{k}$ for the Bratu problem with no preconditioning. Here we fix $\theta =0.9$ and  cap the number of C2D iterations to 4. Form left to right we take $\hat{m} = 2,4,8$. The label refers to the unfolding direction of the core of $X_{k}$. \label{fig:bratu_rank_2}}
\end{center}
\end{figure}

In Figure \ref{fig:bratu_rank_final}, for all six parameter combinations, we plot the singular value distribution for the unfolding of the core in all directions and for the final iterate together with the singular values for the final (and more accurate) solution obtained with the preconditioned Tucker-AA. As can be seen in the figure $\theta =0.9$ has less ``superfluous ranks'' with the tightest fit to the preconditioned solution occurring for C2D capped at 4 iterations. 

\begin{figure}[]
\begin{center}
\includegraphics[width=0.16\textwidth,trim={0.0cm 0.0cm 0.0cm 0.0cm},clip]{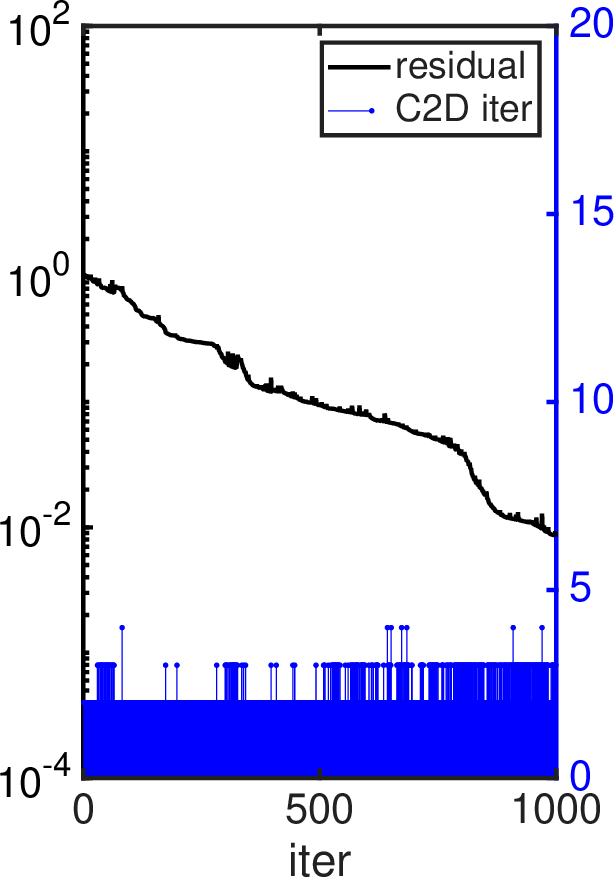} 
\includegraphics[width=0.16\textwidth,trim={0.0cm 0.0cm 0.0cm 0.0cm},clip]{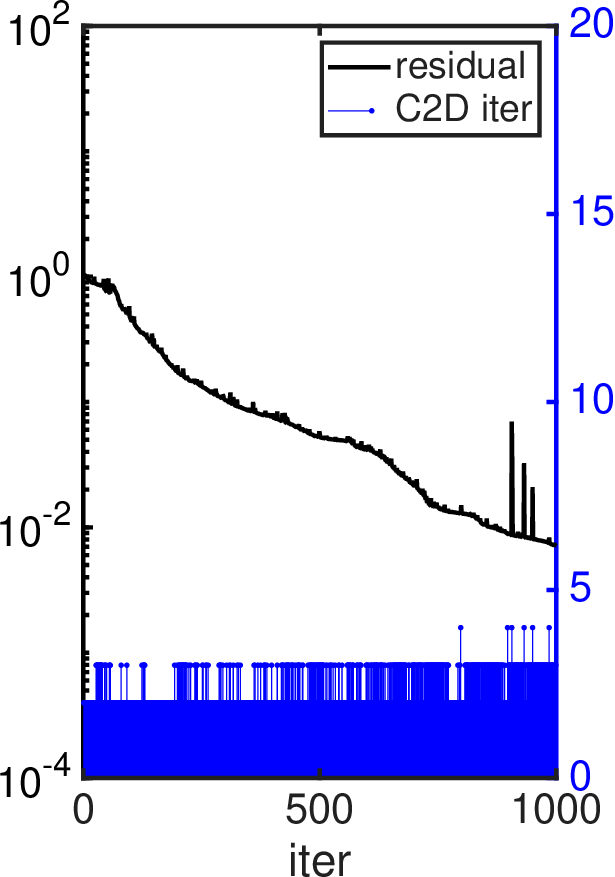} 
\includegraphics[width=0.16\textwidth,trim={0.0cm 0.0cm 0.0cm 0.0cm},clip]{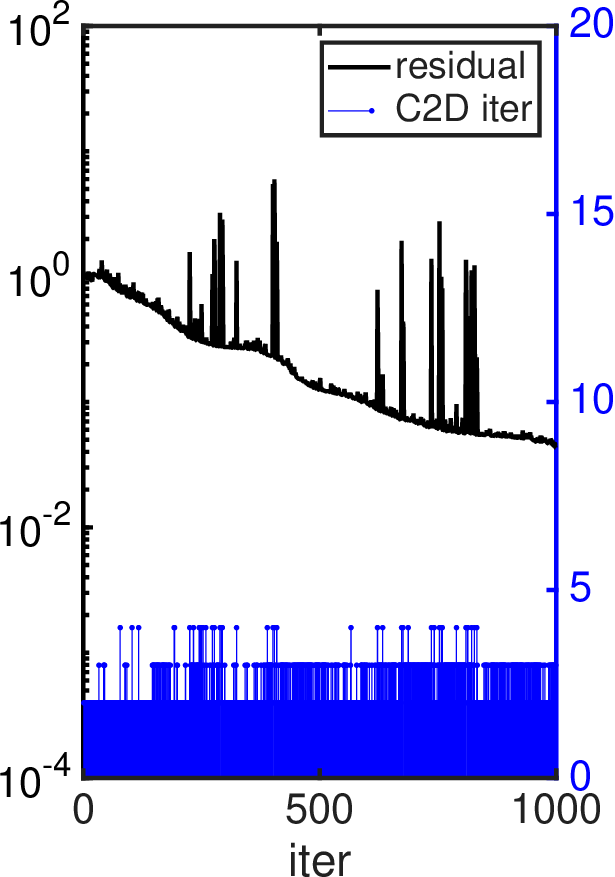} 
\caption{The figure displays the residuals (left vertical axis) and the number of C2D iterations (right axis) of the iterates $X_{k}$ for the Bratu problem with no preconditioning. Here we fix $\theta =0.9$ and  cap the number of C2D iterations to 4. Form left to right we take $\hat{m} = 2,4,8$. \label{fig:bratu_resid_2}}
\end{center}
\end{figure}

\begin{figure}[]
\begin{center}
\includegraphics[width=0.32\textwidth,trim={0.0cm 0.0cm 0.0cm 0.0cm},clip]{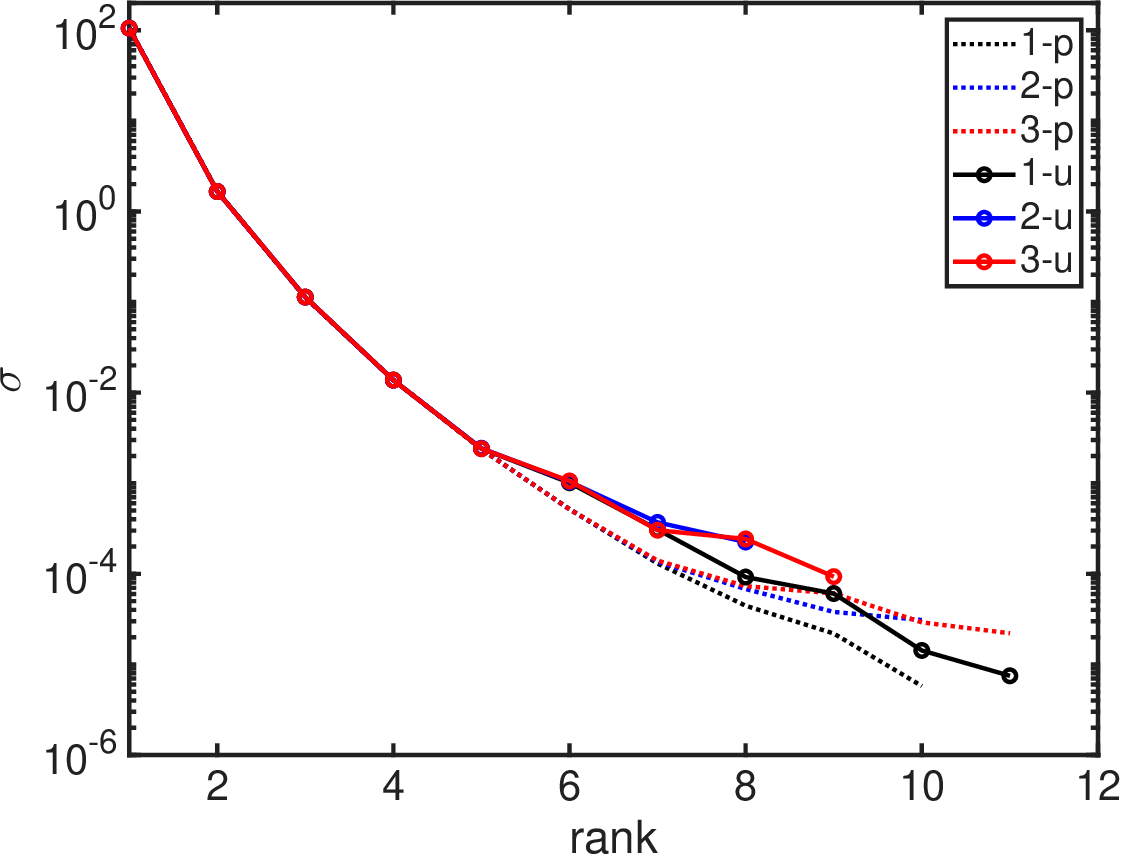} 
\includegraphics[width=0.32\textwidth,trim={0.0cm 0.0cm 0.0cm 0.0cm},clip]{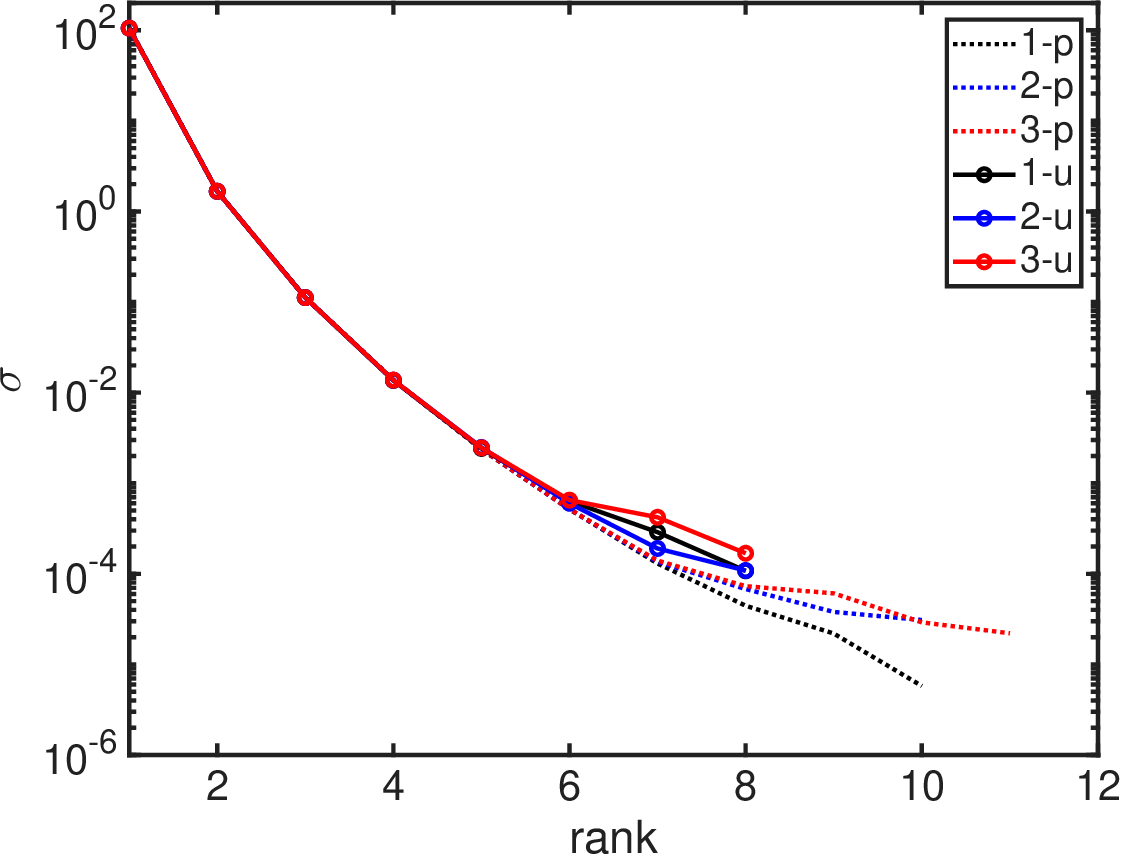} 
\includegraphics[width=0.32\textwidth,trim={0.0cm 0.0cm 0.0cm 0.0cm},clip]{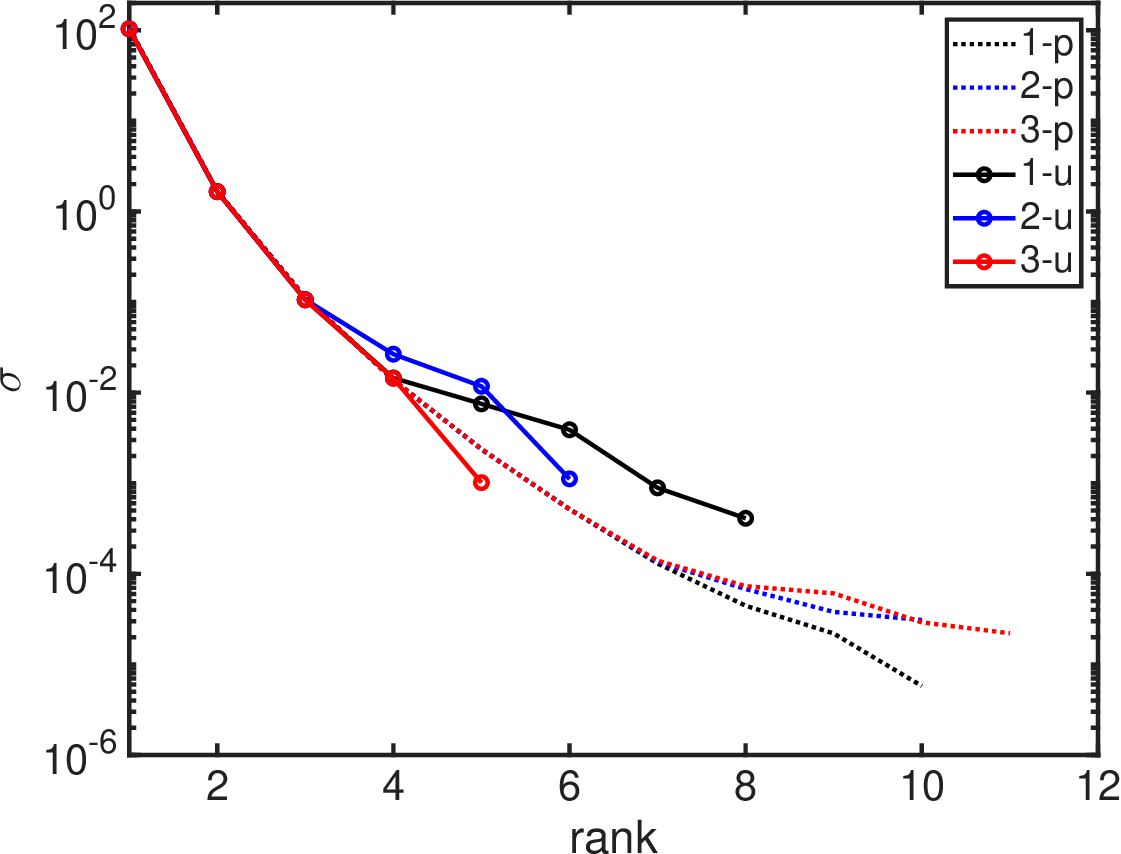} 
\caption{Displayed are the singular value distribution for the unfolding of the core in all directions and for the final iterate together with the singular values for the final   solution obtained with the un-conditioned Tucker-AA. (The dots are obtained with the preconditioned Tucker-AA as reference solution.) The ordering of the graphs is the same as in Figure \ref{fig:bratu_rank_2}. \label{fig:bratu_rank_final_2}}
\end{center}
\end{figure}

For $\theta = 0.9$ and  with C2D capped at 4 iterations we then use $\hat{m} = 2,$ 4 and 8 to check the impact of different window sizes. The results of the same quantities as just discussed are displayed in Figures \ref{fig:bratu_rank_2}, \ref{fig:bratu_resid_2} and \ref{fig:bratu_rank_final_2}. For this test the cases  $\hat{m} = 2$ and  $\hat{m} = 4$ perform similarly for the rank as a function of iteration and for the decrease in residual. The case  $\hat{m} = 4$ appears to have a slightly better distribution of singular values at the final iteration. The case $\hat{m} = 8$ is clearly worse for all the diagnostics. These results align well with our previous experience from the matrix version of Tucker-AA, it is best to keep the window size fairly small.

\subsubsection{The Allen-Cahn equation}
Here we solve the Allen-Cahn equation 
\[
v_t = \nu \Delta v +v - v^3,
\]    
with $\nu=0.01$ and homogenous Neumann boundary conditions in the cube  $(x_1,x_2,x_3) \in [0,2\pi]^3$. The initial data is taken to be   
\[
    v(x_1,x_2,x_3) = \sin(x_1) \sin(x_2) \sin(x_3).
\]
As before $X(i_1,i_2,i_3) \approx v(i_1h, i_2h ,i_3h)$ on an equidistant grid with gridsize $h = 1/(n+1)$ in all directions. We adopt the solver from Section \ref{sec:PoissonDFT} as a preconditioner (since we have Neumann boundary conditions we instead use the discrete cosine transform). We discretize in time using the backward Euler method, resulting in a scheme defined by 
\begin{gather}
\begin{split}
H_{\rm AC}(i_1,i_2,i_3;X) &= X(i_1+1,i_2,i_3) - X_{old}(i_1+1,i_2,i_3) - \Delta t \Big(  \\
& + \frac{\nu }{h^2}\left(X(i_1+1,i_2,i_3) - 2X(i_1,i_2,i_3) + X(i_1-1,i_2,i_3)  \right) \\
&+ \frac{\nu }{h^2} \left(X(i_1,i_2+1,i_3) - 2X(i_1,i_2,i_3) + X(i_1,i_2-1,i_3)  \right) \\
&+\frac{\nu }{h^2} \left(X(i_1,i_2,i_3+1) - 2X(i_1,i_2,i_3) + X(i_1,i_2,i_3-1)  \right) \\
&+ X(i_1,i_2,i_3) - X(i_1,i_2,i_3)^3 \Big),
\end{split}
\end{gather}
with appropriate modifications at boundaries. Here $X_{old}$ refers to the solution at the current time. The fixed point function in Tucker-AA, $H(i_1,i_2,i_3)$,  is obtained by applying the preconditioned Richardson iteration. We have 
\[
X_{k+1}(i_1,i_2,i_3) = H(i_1,i_2,i_3;X_{k},\alpha) \equiv  X_{k}(i_1,i_2,i_3) + \alpha M(H_{\rm AC}(i_1,i_2,i_3;X_{k})).
\] 

In each timestep the Tucker-AA parameters used are the following ${\rm TOL}=5\cdot 10^{-2} \rho_0, \hat{m}=3, \theta=0.9$ and $\alpha = 0.2$. The current timestep solution is used as initial guess for the next timestep. We use $n = 255$ and take 50 timesteps to advance to time 5.

\begin{figure}[]
\begin{center}
\includegraphics[width=0.35\textwidth,trim={0.0cm 0.0cm 0.0cm 0.0cm},clip]{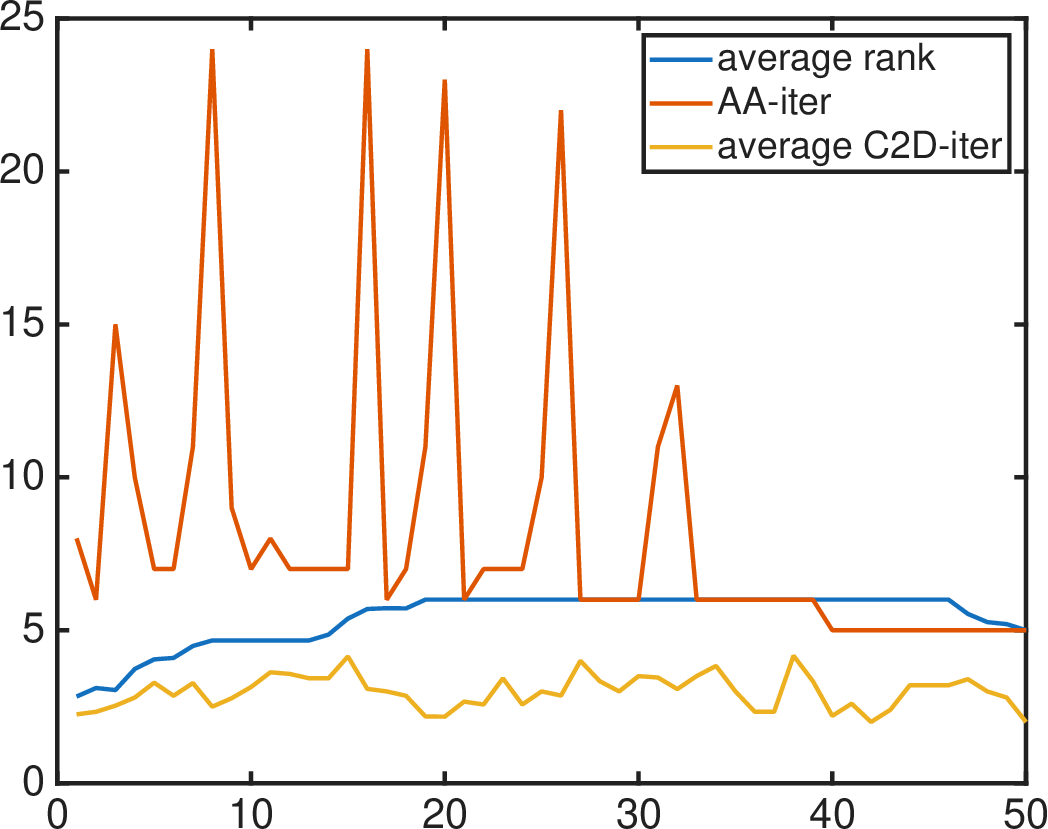} 
\caption{Displayed are the average rank of the solution to Allen-Cahn equation example, the average number of C2D iterations per Tucker-AA solve, and the number of iterations for each Tucker-AA solve.  \label{fig:allen_cahn_iter_vs_timestep}}
\end{center}
\end{figure}

\begin{figure}[]
\begin{center}
\includegraphics[width=0.32\textwidth,trim={3.7cm 0.0cm 2.4cm 0.0cm},clip]{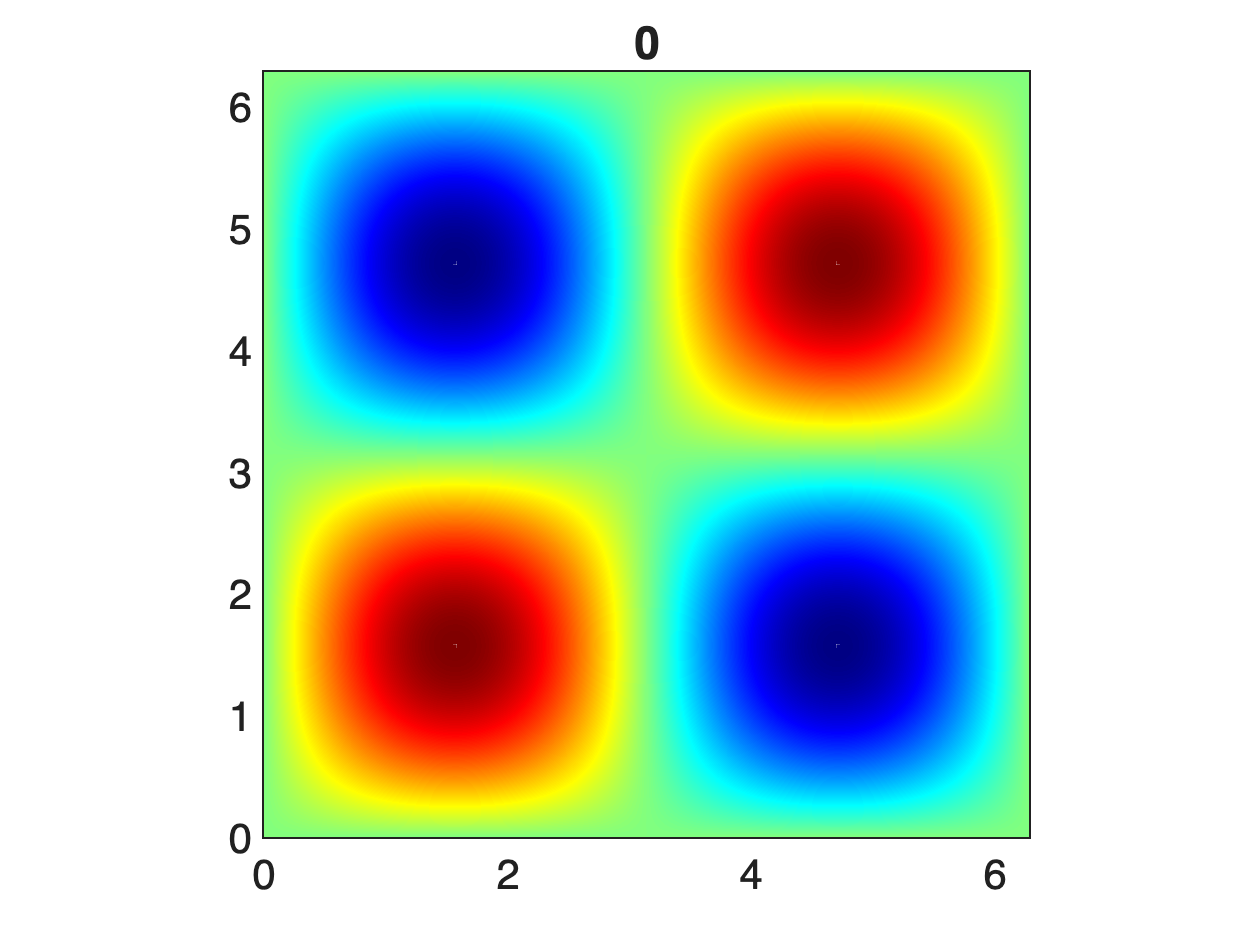} 
\includegraphics[width=0.32\textwidth,trim={3.7cm 0.0cm 2.4cm 0.0cm},clip]{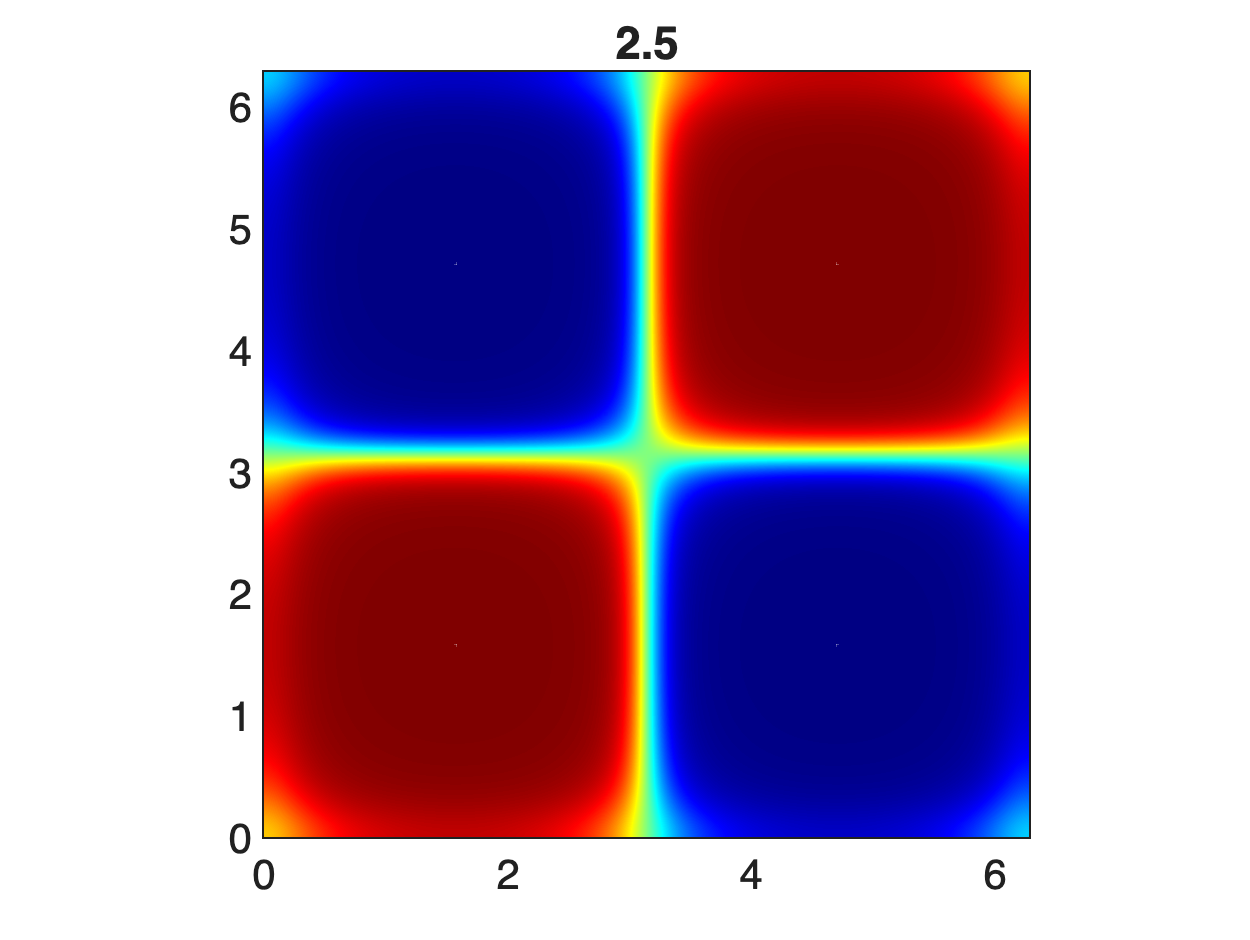} 
\includegraphics[width=0.32\textwidth,trim={3.7cm 0.0cm 2.4cm 0.0cm},clip]{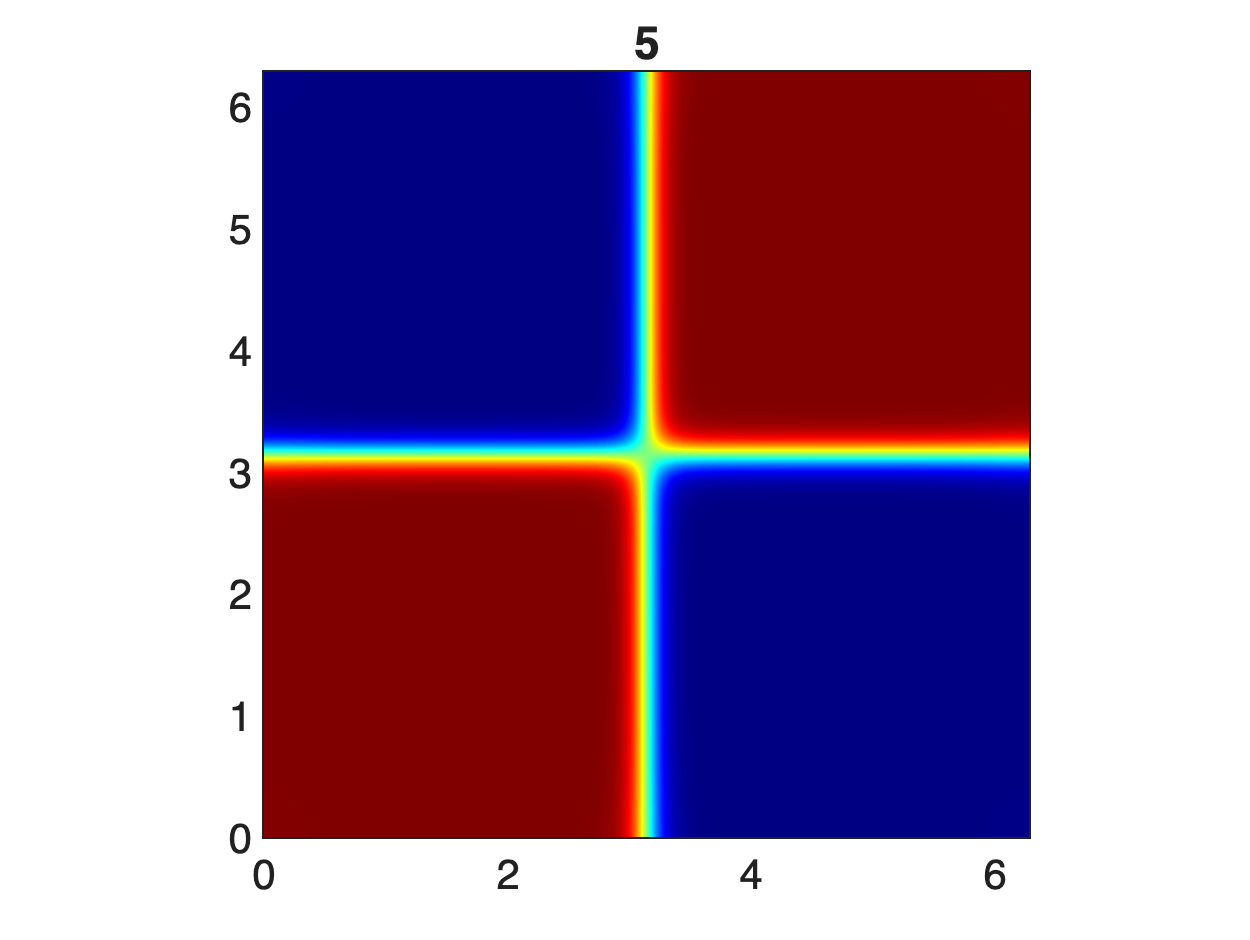} 
\caption{Snapshots of the solution to the Allen-Cahn example at times 0, 2.5 and 5. Here the color scale ranges from -1 (blue) to 1 (red). The figures are at the cut $x_3 = \pi/2$. \label{fig:allen_cahn_sol}}
\end{center}
\end{figure}

In Figure \ref{fig:allen_cahn_iter_vs_timestep} we display the average rank, the number of Tucker-AA iterations and the average number of C2D iterations as a function of the timestep. As can be seen the average rank remains small as does the average number of C2D iterations. The iteration numbers of Tucker-AA are also small but vary more. In Figure \ref{fig:allen_cahn_sol} we display the approximate solution at a cut $x_3 = \pi/2$ at time 0, 2.5 and 5. At time 5 the solution is essentially stationary. 

\section{Conclusions}
\label{sec:conclude}
In this paper, we proposed Cross$^2$-DEIM as a fast computational approach for cross approximation in Tucker tensor format. This method resembles FSTD2 sampling, which samples $O(r)$ fibers in each mode.   $O(dnr+r^d)$ entries are sampled within each iteration count. The computational cost of Cross$^2$-DEIM   besides fiber sampling cost is $O(dn r^2+r^{d+1})$ which is   lower than the existing ones in the literature  \cite{ahmadi2021cross} if the iteration number on $k$ stays $O(1).$ We demonstrate the performance of Cross$^2$-DEIM for tensor approximation and low-rank decompositions to fundamental solution to Helmholtz equation.  We used Cross$^2$-DEIM to devise a sublinear fast direct Poisson solver in high dimensions, which can serve as a preconditioner for low-rank methods for elliptic problems.

We designed a low-rank solver for nonlinear tensor equation in Tucker format: Tucker-AA. It is based on the AA framework with warm-started Cross$^2$-DEIM and scheduling of truncation parameters. The Tucker-AA algorithm shows good performance  for applications to numerical discretizations of nonlinear PDEs.

Future work includes considering tensor train and various other tensor formats, and refinement and analysis of the proposed schemes and their applications.

\bibliography{lraa}
\bibliographystyle{plain}
\end{document}